\documentclass{amsart}
\newtheorem{lemma}{Lemma}[section]
\newtheorem{propos}[lemma]{Proposition}

\newtheorem{theorem}[lemma]{Theorem}
\newtheorem{cor}[lemma]{Corollary}

\newtheorem{defin}[lemma]{Definition}
\theoremstyle{remark}

\newcommand{\cg}{\mathfrak{g}}

\newcommand{\C}{\mathbb{C}}

\newcommand{\del}{\partial}

\newcommand{\extd}{{\rm d}}
\newcommand{\isom}{{\cong}}
\newcommand{\eps}{{\epsilon}}
\newcommand{\tens}{\mathop{\otimes}}
\newcommand{\la}{{\triangleright}}
\newcommand{\ra}{{\triangleleft}}

\newcommand{\Ad}{{\rm Ad}}
\newcommand{\ad}{{\rm ad}}
\newcommand{\id}{{\rm id}}
\newcommand{\<}{\langle}
\renewcommand{\>}{\rangle}

\renewcommand{\o}{{}_{\scriptscriptstyle(1)}}
\renewcommand{\t}{{}_{\scriptscriptstyle(2)}}
\newcommand{\thr}{{}_{\scriptscriptstyle(3)}}
\newcommand{\fo}{{}_{\scriptscriptstyle(4)}}
\newcommand{\bz}{{}_{{\scriptscriptstyle(0)}}}
\newcommand{\bt}{{}_{{\scriptscriptstyle(\infty)}}}
\newcommand{\bo}{{}_{{\scriptscriptstyle(1)}}}
\newcommand{\so}{{}^{\scriptscriptstyle[1]}}
\newcommand{\st}{{}^{\scriptscriptstyle[2]}}
\newcommand{\sth}{{}^{\scriptscriptstyle[3]}}

\newcommand{\lbiprod}{{>\!\!\!\triangleleft\kern-.33em\cdot}}
\newcommand{\rbiprod}{{\cdot\kern-.33em\triangleright\!\!\!<}}

\newcommand{\eproof}{$\quad \square$\bigskip}
\newcommand{\nquad}{\kern-60pt}

\newcommand{\eqn}[2]{\begin{equation}#2\label{#1}\end{equation}}

\begin{document}
\title[Semi-classical differential structures]{\rm\large
SEMI-CLASSICAL DIFFERENTIAL STRUCTURES}

\author{E.J. Beggs + S. Majid}%
\address{EJB: Department of Mathematics, University of Wales
   Swansea, SA2 8PP, UK\\
SM: School of Mathematical Sciences\\
Queen Mary, University of London\\ 327 Mile End Rd,  London E1
4NS, UK}

\thanks{SM is a Royal Society University Research Fellow}%


\maketitle

\begin{abstract} We semiclassicalise the standard notion of differential
calculus in noncommutative geometry on algebras and quantum
groups. We show in the symplectic case that the infinitesimal data
for a differential calculus is a symplectic connection, and
interpret its curvature as lowest order nonassociativity of the
exterior algebra. Semiclassicalisation of the noncommutative torus
provides an example with zero curvature. In the Poisson-Lie group
case we study left-covariant infinitesimal data in terms of
partially defined preconnections. We show that the moduli space of bicovariant
infinitesimal data for quasitriangular Poisson-Lie groups has a
canonical reference point which is flat in the triangular case.
Using a theorem of Kostant, we completely determine the moduli
space when the Lie algebra is simple: the canonical preconnection is the unique point for other than $sl_n$, $n>2$, when
the moduli space is 1-dimensional. We relate the canonical preconnection to Drinfeld twists and thereby quantise it to a super
coquasi-Hopf exterior algebra. We also discuss links with Fedosov
quantisation.
\end{abstract}

\section{Introduction}

Usually the quantisation problem in physics consists of a
commutative algebra of functions equipped with a Poisson bracket
and the search for a noncommutative algebra with commutators
reproducing this to lowest order in a deformation parameter
$\hbar$. It is well known that actually the converse problem is
more well posed: given a noncommutative algebra which is a flat
deformation one may semiclassicalise its structure and recover the
Poisson bracket of which it is a quantisation. Either way Poisson
brackets are the semiclassical data for associative noncommutative
algebras.

In this second point of view the starting point is the
noncommutative algebra itself. In the last two decades the
`geometry' of such algebras has been well developed in different
approaches, such as from operator theory\cite{Con} or Hopf
algebras\cite{Ma:book}, and in particular the notion of
differential structures, quantum bundles, Riemannian structures
etc have been fairly well established from an intrinsically
noncommutative algebraic point of view. In this paper and its
sequel we will semiclassicalise these various notions from
noncommutative geometry to elucidate the classical infinitesimal
data of which they are quantisations. In this first part we limit
ourselves to the differential calculus and differential forms. In
a sequel we will proceed to quantum bundles and Riemannian
structures at the semiclassical level.

Let $A$ be an algebra. The by-now standard notion of differential
calculus in noncommutative geometry is to specify an
$A-A$-bimodule $\Omega^1$ of `1-forms' and a linear map
$\extd:A\to \Omega^1$ obeying:

\begin{enumerate}
\item Leibniz rule $\extd(ab)=a(\extd b)+ (\extd a)b$ for all
$a,b\in A$
\item Surjectivity, i.e. $\{a\extd b\}$ span $\Omega^1$
\end{enumerate}

We do not demand that $a \extd b= (\extd b)a$ i.e. that 1-forms
and functions commute, just as we do not demand that $A$ is
commutative. The above is a first order calculus and we will find
that its semiclassical data is a compatible partially defined `preconnection' (a
symplectic torsion free connection in the symplectic case when the
symplectic form is central to lowest order), see
Propositions~\ref{peqq2} and~\ref{omegacen}.

We can ask further for an entire differential graded `algebra' of
forms of all degree and $\extd $ such that $\extd ^2=0$. Usually
one demands an associative such exterior algebra and this case corresponds to a flat symplectic connection or preconnection. The
formal analysis for both results is in Section~2 and the geometric
meaning is in Section~3. We find in particular three super Jacobi
identity obstructions $J_1,J_2,J_3$. The geometric meaning of the
$J_1$ obstruction turns out to be the curvature of the preconnection. 
We give geometric conditions for the others also
(involving now the torsion). An example of a flat connection is
provided by the noncommutative torus at the semiclassical level.
\footnote{Note added in press:   in fact the requirement of  flatness in the usual associative case has been found earlier in the paper \cite{Eli}  using a different notion of `contravariant connections' (this does not affect our main results which are intended to reach beyond the flat case). We thank E. Hawkins for bringing this work to our attention and for helpful comments.}

Next, in Section~4, we specialise to the case where $A$ is a Hopf
algebra or `quantum group'. These provide examples of
noncommutative geometry which are well controlled though the
requirements of the `group' structure. Their semiclassicalisation
was worked out in the 1980's by V.G. Drinfeld as the notion of a
Poisson-Lie group. In this context it is natural to restrict to
differential structures that are left, right or
bicovariant\cite{Wor:dif}. We analyse the left covariance
restrictions on a Poisson-Lie group $G$ at the semiclassical level
in Sections~4.1 and~4.2 and give a formulation of the data in Lie
algebraic terms of a map $\Xi:\cg\to\cg\tens\cg$. We give an
example on the group manifold $SU_2$ in Section~4.3. After that,
Secton~4.4 analyzes the stronger requirement of a bicovariant
calculus at the semiclassical level, where Proposition~4.15
together with Corollary~4.7 reduce the classification to
$\Ad$-invariant symmetric maps $\hat\Xi:\cg\to {\rm Sym}^2(\cg)$
when the Poisson-Lie group is quasitriangular (this includes the
standard $q$-deformation quantum groups at the semiclassical
level). We conclude in this case (Theorem~4.18 and
Proposition~4.19) that there is a canonical choice and it has
curvature given by the Schouten bracket $[[r_-,r_-]]$, where $r_-$
is the antisymmetric part of the quasitriangular structure $r\in
\cg\tens \cg$. Thus the classical Yang-Baxter equation has a
direct geometric meaning as curvature. We show (Theorem~4.20) that
for all simple $\cg$ other than the $sl_n$ series, $n>2$, it is
the {\em only} choice for a first order bicovariant semiclassical
calculus (i.e. a compatible preconnection). Moreover, our
results imply that this extends associatively to lowest order if
and only if $\cg$ is triangular. This fits with the `quantum Lie
functor' for a canonical differential calculus on quasitriangular
Hopf algebras introduced in \cite{GomMa} which is known to
typically give trivial answers unless the quantum group is
triangular. Our results prove indeed that a strictly
quasitriangular Poisson-Lie group has {\em no} bicovariant
calculus that is a deformation (with the same dimensions etc.) of
the classical calculus in a high degree of generality. This
explains the typical experience in quantum group theory
 where strictly quasitriangular quantum groups
typically do not admit bicovariant differential calculi of the
usual classical dimension but require some form of central or
other extension. What is remarkable is that this obstruction
identified in our analysis is not at the semiclassical first order
level where we find the perfectly good canonical 
preconnection above; it enters at the semiclassical super Jacobi
identity level concerning the associativity in the bimodule and
exterior algebra structures.

This last point is taken up in Section~5, where we provide a
different point of view on the canonical preconnection for
the standard quasitriangular Poisson-Lie groups, now inspired by
Drinfeld's theory of quasi-Hopf algebras\cite{Dri}. Drinfeld
showed (effectively) that the standard quantum groups $C_\hbar(G)$
in a suitable deformation setting are isomorphic to the twisting
by a cochain of a coquasi-Hopf algebra structure on the classical
group function algebras. Even though the cochain is not a cocycle
it happens that the resulting algebra is associative. This
accident does not extend to the exterior algebra, i.e. when the
classical exterior algebra is similarly twisted as in
\cite{MaOec:twi} it becomes nonassociative. We semiclassicalise
this construction and understand our canonical preconnection
in these terms as corresponding to a quasiassociative calculus
(Proposition~5.3). Conversely put, we quantise the canonical
preconnection to a super coquasi-Hopf exterior algebra
$\Omega(C_\hbar(G))$ by these means.

Finally, we note that our main result that the data for a
semiclassical calculus is a  pair consisting of a Poisson bracket
and a compatible connection or preconnection is, in the symplectic case,
precisely the data used by Fedosov \cite{Fed} to solve the
quantisation problem. We make some remarks about this in
Section~6. Our result suggests that in the Fedosov construction
there is also for free a noncommutative differential calculus also
constructed from the input data and with semiclassical limit the
Fedosov symplectic connection. If so then whereas the problem of
quantising a symplectic structure to an algebra is not unique as
one must chose the connection, quantising the pair consisting of a
symplectic structure and connection to an algebra and differential
calculus would be unique. Also, because of the role played in
Hamiltonian mechanics by symplectic structures, we expect a
physical role also of the symplectic connection and here its
interpretation as controlling the quantisation of the differential
structure suggests a possible role.

In this paper we use the following conventions for the
  torsion and curvature tensors.
For vector fields $U$, $V$ and $W$ we have
\[
T(U,V)^k = T^k_{ij}\,U^i\,V^j\ ,\quad
T(U,V)\,=\,\nabla_U\,V-\nabla_V\,U\,-\,[U,V]\ \]\[ R(U,V)(W)^l =
R^l_{\phantom{l}ijk}\, W^i\, U^j\,V^k\ ,\quad R(U,V)(W)
\,=\,\nabla_U\nabla_V\,W-\nabla_V\nabla_U\,W-\nabla_{[U,V]}\,W\ .
\]
In terms of Christoffel symbols we have
\[
(\nabla_V\,W)^l\,=\,V^k\,(W^l_{\phantom{l},k}+\Gamma^l_{ki}\,W^i)\
,\quad T^k_{ij}\,=\,\Gamma^k_{ij}-\Gamma^k_{ji}\ \]
\[
R^{l}_{\phantom{l}ijk}\,=\, \frac{\del \Gamma^l_{ki}}{\del
x^j}\,-\, \frac{\del \Gamma^l_{ji}}{\del x^k}\,+\,
\Gamma^m_{ki}\,\Gamma^l_{jm}\,-\,\Gamma^m_{ji}\,\Gamma^l_{km}\ .
\]

\section{Deformation analysis of noncommutative differential structures}

In this section we perform the required algebraic deformation
analysis and prove some basic lemmas. Their geometric meaning will
then be explored in later sections.

As usual, we deform commutative multiplication on $C^\infty(M)$ on
a classical manifold $M$ to the associative multiplication
$x\bullet y$, where $x\bullet y\,=\,x\,y+O(\hbar)$. If we assume
that the commutator can be written as $[x,y]_\bullet=x\bullet y-
y\bullet x=\,\hbar\,\{x,y\}+O(\hbar^2)$, we see that
$\{x,y\}=-\{y,x\}$ and (by considering the two ways of writing
$z\bullet x\bullet y- x\bullet y\bullet z$ to first order in
$\hbar$), $\{z,xy\}\,=\, \{z,x\}\,y\,+\,x\{z,y\}$, i.e.\, $\{\ ,\
\}$ is a Poisson bracket. Formally speaking, the deformed algebra
$C^\infty_\hbar(M)$ can be formulated as an algebra over
$\C[[\hbar]]$, topologically free with $C^\infty_\hbar(M)/\hbar
C^\infty_\hbar(M)\isom C^\infty(M)$. However, for most of the paper we
actually need only that structure maps depend on an $\hbar$  parameter
permitting of power series expansion and comparison of lower order
terms as specified, which is therefore the line that we take. 
The same remark applies for differential forms
in what follows.

\subsection{Deformation of the classical bimodule structure}
In this subsection we begin with assumptions on the associativity
of the differential calculus, and see what this means in terms of
the super-commutators $[\, ,]_\bullet$. We define a {\em formal
deformation} of the differential calculus on $M$ to lowest order
in $\hbar$ to be the following:

As sets the functions $C^\infty(M)$ and $n$-forms $\Omega^n(M)$
take their classical values.  The symbol $\bullet$ is used for the
deformed multiplication. $\Omega^1(M)$ is a bimodule over
    $(C^\infty(M),\bullet)$ to order $O(\hbar^2)$, i.e.\ for $x,y\in
    C^\infty(M)$ and $\tau\in\Omega^1(M)$
\begin{eqnarray}\label{appbimod}
&&(x\bullet y)\bullet \tau\,-\,x\bullet(y\bullet
\tau)\,=\,O(\hbar^2)\nonumber \\
&&\tau\bullet(x\bullet y)\,-\,(\tau\bullet x)\bullet
y\,=\,O(\hbar^2)\nonumber\\ &&(x\bullet \tau)\bullet
y\,-\,x\bullet(\tau\bullet y) = O(\hbar^2)\ . \end{eqnarray}

We also suppose that $x\bullet\tau\,=\,x\,\tau\,+\,O(\hbar)$ and
$\tau\bullet x\,=\,x\,\tau\,+\,O(\hbar)$, and hence define
$\gamma$ by \eqn{gamma}{x\bullet\tau-\tau\bullet
x=[x,\tau]_\bullet=\hbar\,\gamma(x,\tau)+O(\hbar^2).}

We make the assumption that the deformed exterior multiplication
$\wedge_\bullet$ is associative to $O(\hbar^{2})$. Note that where
one of the forms is a zero form, we just use $\bullet$ instead of
$\wedge_\bullet$.  We also assume that
$\tau\wedge_{\bullet}\eta=\tau\wedge\eta+O(\hbar)$.

The deformed $\extd $ operator $\extd
^\bullet:\Omega^n(M)\to\Omega^{n+1}(M) $, is related to the usual
$\extd $ by $\extd ^\bullet x=\extd x+O(\hbar)$.  Also $\extd
^\bullet$ is a graded derivation to order $O(\hbar^2)$, i.e.\
\begin{eqnarray}\label{appder}
\extd ^\bullet(\xi\wedge_{\bullet} \eta)\,=\,\extd ^\bullet
\xi\wedge_{\bullet} \eta\,+\,(-1)^{{\rm deg}\xi}\,
\xi\wedge_{\bullet} \extd ^\bullet \eta\,+\,O(\hbar^2)\ .
\end{eqnarray}

\begin{propos}\label{peqq2}
The map $\gamma(-,\tau)$ is a derivation on $C^\infty(M)$ for all
$\tau\in\Omega^1(M)$ and obeys \eqn{eqq1}{\gamma(x,\tau\, y)\,=\,
\gamma(x,\tau)\, y\,+\, \tau\,\{x,y\}\ .} We call any  map
$\gamma$ with these properties a  {\em preconnection}.
Moreover, $\gamma$ is compatible with the Poisson structure in the
sense \eqn{eqq2}{\extd \{x,y\}\,=\,\gamma(x,\extd
y)\,-\,\gamma(y,\extd x)\ .}
\end{propos} \proof From the approximate bimodule rules
(\ref{appbimod}) we have, for $x,y\in C^\infty(M)$ and
$\tau\in\Omega^1(M)$.
\begin{eqnarray*}
[x\bullet y,\tau]_\bullet\,=\,x\bullet[y,\tau]_\bullet\,+\,
[x,\tau]_\bullet \bullet y\,+\,O(\hbar^2)\ , \end{eqnarray*} from
which we deduce
\begin{eqnarray}\label{eqq3}
\gamma(xy,\tau)\,=\,x\,\gamma(y,\tau)\,  \,+\, \gamma(x,\tau)\,y\
. \end{eqnarray} The next formula is deduced from
\begin{eqnarray*}
[x,\tau\bullet y]_\bullet\,=\, [x,\tau]_\bullet\bullet
y\,+\,\tau\bullet[x,y]_\bullet\,+\,O(\hbar^2)\ , \end{eqnarray*}
and the last formula from the approximate derivation rule
(\ref{appder})
\begin{eqnarray*}
\extd ^\bullet[x,y]_\bullet\,=\,[\extd ^\bullet x,y]_\bullet\,+\,
[x,\extd ^\bullet y]_\bullet\,+\,O(\hbar^2)\ .\quad\square
\end{eqnarray*}

The terminology for $\gamma$ is justified in Section~3 below.

\begin{propos} The commutator of a function $x$ with higher forms is
given by the following formula, which can also be viewed as the
extension of $\gamma$ to higher forms:
\[
[x,\tau\wedge_{\bullet}\eta]_{\bullet}\,=\,\hbar\,\gamma(x,\tau\wedge\eta)
\,+\, O(\hbar^{2})\,=\,\hbar\,(\gamma(x,\tau)\wedge\eta
+\tau\wedge\gamma(x,\eta))\,+\, O(\hbar^{2})\ .
\]
\end{propos} \proof This can be seen by rearranging the formula
\[
[x,\tau\wedge_{\bullet}\eta]_{\bullet}\,=\,
x\bullet(\tau\wedge_{\bullet}\eta)\,-\,(\tau\wedge_{\bullet}\eta)\bullet
x \,+\,O(\hbar^{2})\ .  \square
\]

\subsection{Super commutator of forms and Jacobi terms}
If the operation $\wedge_\bullet$ on $\Omega^*$ is associative (up
to $O(\hbar^2)$) and super-commutative (up to $O(\hbar)$), then
the super-Jacobi identities for the super-commutator $[,]_\bullet$
hold up to $O(\hbar^2)$. In this subsection we ask what conditions
would be necessary for the super-Jacobi identities to hold up to
$O(\hbar^3)$.

\begin{defin} Define $J_{i}:C^{\infty}(M)^{\tens 3}\to \Omega^{i}(M)$
     by \begin{eqnarray*}
[x,[y,\extd z]_\bullet ]_\bullet \,+\,[\extd z,[x,y]_\bullet
]_\bullet \,+\,[y,[\extd z,x]_\bullet ] _\bullet &=& \hbar^2\
J_1(x,y,z)\,+\,O(\hbar^3)\ ,\cr [x,[\extd y,\extd z]_\bullet
]_\bullet \,-\,[\extd z,[x,\extd y]_\bullet ]_\bullet \,+\,[\extd
y,[\extd z,x]_\bullet ] _\bullet &=& \hbar^2\
J_2(x,y,z)\,+\,O(\hbar^3)\ ,\cr [\extd x,[\extd y,\extd z]_\bullet
]_\bullet \,+\,[\extd z,[\extd x,\extd y]_\bullet ]_\bullet
\,+\,[\extd y,[\extd z,\extd x]_\bullet ] _\bullet &=& \hbar^2\
J_3(x,y,z)\,+\,O(\hbar^3)\ . \end{eqnarray*} \end{defin}

\begin{propos} \label{oouu}

1)\quad If $\extd \,J_{2}$ vanishes identically then so does
$J_{3}$.

2)\quad If $\extd \,J_{1}$ vanishes identically then $J_{2}$ is
totally symmetric in its 3 arguments.
\end{propos}
\proof By applying $\extd ^\bullet$ we see that $\extd
\,J_2(x,y,z)=J_3(x,y,z)+O(\hbar)$, so the vanishing of $\extd
\,J_2$ implies the vanishing of $J_3$.  For the second statement,
note that $\extd J_1(x,y,z)=J_2(x,y,z)-J_2(y,z,x)+O(\hbar)$, and
by combining this with the more obvious identity
$J_2(x,z,y)=J_2(x,y,z)$ we find that the vanishing of $\extd J_1$
implies that $J_2(x,y,z)$ is totally symmetric in $x$, $y$ and
$z$. \eproof

\begin{propos}\label{jac2} Suppose that the first super Jacobi identity (for
two functions and a 1-form) holds to $O(\hbar^2)$, that
$J_{2}(x,x,x)=0$ for all $x\in C^\infty(M)$, and that the
following conditions are satisfied for all $x,y\in C^{\infty}(M)$,
$\tau,\eta\in\Omega^{1}$ and $\pi\in\Omega^{2}$:
\begin{eqnarray*}
     [\tau,x\,\eta]_{\bullet} &=& [\tau,x]_{\bullet}\wedge \eta
     \,+\, x\,[\tau,\eta]_{\bullet}\,+\,O(\hbar^{2})\ ,\cr
     [y,x\,\pi]_{\bullet} &=& [y,x]_{\bullet} \,\pi\,+\,
     x\,[y,\pi]_{\bullet}\,+\,O(\hbar^{2}) \ ,\cr
     [y,\tau\wedge\eta]_{\bullet} &=& [y,\tau]_{\bullet}\wedge\eta
     \,+\, \tau\wedge [y,\eta]_{\bullet}\,+\,O(\hbar^{2})\ .
\end{eqnarray*}
Then the second super Jacobi identity (for a function and two
1-forms) holds to $O(\hbar^2)$.
\end{propos} \proof As the super Jacobi identity for two functions
and a 1-form holds to $O(\hbar^2)$,
  we see that $J_1$ vanishes identically. By Proposition~\ref{oouu} $J_2$ is
  completely symmetric, and as $J_{2}(x,x,x)=0$ for all $x\in C^\infty(M)$ we
deduce that $J_2$ is identically zero. Now use the fact that
linear combinations of the form $x\,\extd y$ span $\Omega^{1}$.
\eproof

\begin{propos}\label{jac3} Suppose that the conditions for Proposition
\ref{jac2} hold, and that the following conditions hold for all
$x\in C^{\infty}(M)$, $\tau,\eta,\xi\in\Omega^{1}$ and
$\pi\in\Omega^{2}$:
  \begin{eqnarray*}
     [\tau,\xi\wedge\eta]_{\bullet} &=& [\tau,\xi]_{\bullet}\wedge
  \eta \,-\,
     \xi\wedge[\tau,\eta]_{\bullet}\,+\,O(\hbar^{2})\ ,\cr
     [\tau,x\,\pi]_{\bullet} &=& [\tau,x]_{\bullet} \wedge\pi\,+\,
     \tau\wedge[y,\pi]_{\bullet}\,+\,O(\hbar^{2}) \ ,\cr
     [x\,\tau,\pi]_{\bullet} &=& [x,\pi]_{\bullet} \wedge \tau\,+\,
     x\,[\tau,\pi]_{\bullet}\,+\,O(\hbar^{2}) \ .
\end{eqnarray*}
  Then the third
super Jacobi identity (for three 1-forms) holds to $O(\hbar^2)$.
\end{propos} \proof Use the fact that linear
combinations of the form $x\,dy$ are dense in $\Omega^{1}$.
\eproof

\section{Geometric interpretation of the semiclassical data}

Here we look at the geometric meaning of the map $\gamma$ in the
semiclassical data. The full picture emerges in the symplectic
case as a symplectic connection, but first some remarks about the
general case. For general Poisson bracket we have only a
preconnection or `partially defined connection' which we denote $\hat\nabla_{x}$. It
should be thought of and sometimes is a usual covariant derivative $\nabla_{\hat x}$ along Hamiltonian vector fields $\hat x=\{x,-\}$
associated to $x\in C^\infty(M)$, and is defined by \eqn{nablagamma}{
\hat\nabla_{ x}\tau =\gamma(x,\tau),\quad \forall
\tau\in\Omega^1(M).} Indeed, the derivation property of a
connection on Hamiltonian vector fields is
\[ \hat\nabla_{ x}(y\tau)=y\hat\nabla_{
x}\tau+\hat x(y)\tau\] which is (\ref{eqq1}) in
Proposition~\ref{peqq2}. Meanwhile, writing
$\{xy,z\}=y\{x,z\}+x\{y,z\}$ for all $z$ as $\widehat{xy}=y\hat
x+x\hat y$, the derivation property of $\gamma(-,\tau)$ plays the
role of the usual tensoriality of a connection with respect to
vector field direction of differentiation.

Similarly, we define the curvature of a preconnection in the
usual way but only on such vector fields, where \eqn{partialR}{
R(\hat x,\hat y)=\hat\nabla_{ x}\hat\nabla_{ y}-
\hat\nabla_{y}\hat\nabla_{ x}-\hat\nabla_{\{x,y\}}} given that the Jacobi
identity for the Poisson bracket means that $[\hat x,\hat
y]=\widehat{\{x,y\}}$. Its deformation-theoretic meaning is:

\begin{propos}\label{curv} For a semiclassical differential calculus,
\[ [x,[y,\tau]_\bullet]_\bullet\,+\,
[\tau,[x,y]_\bullet]_\bullet\,+\, [y,[\tau,x]_\bullet]_\bullet
=\hbar^2 R(\hat x,\hat y)(\tau)+O(\hbar^3)\ ,
\] i.e.\ the obstruction
to the first super Jacobi identity is $J_1(x,y,z)=R(\hat x,\hat
y)(\extd z)$.
\end{propos}
\proof  First we calculate
\[
[x,[y,\tau]_\bullet]_\bullet\,=\,[x,\hbar\,\hat\nabla_{
y}\tau+O(\hbar^2)]_\bullet \,=\, \hbar^2\,\hat\nabla_{
x}\hat\nabla_{ y}\tau+O(\hbar^3)\ .
\]
Applying this into the first super Jacobi identity with two
functions $x,y$ and a 1-form $\tau$ and the definition of $\gamma$
gives the result. \eproof

In the same way, we define the torsion tensor again as usual but
partially, by \eqn{partialT}{ T(\hat x,\hat y)=\hat\nabla_{ x}\hat
y-\hat\nabla_{ y}\hat x-\widehat{\{x,y\}}} where $\hat\nabla_{ x}$
on vector fields is defined as usual via
\[ \<\hat\nabla_{ x}v,\tau\>=\hat x(\<v,\tau\>)-\<v,\hat\nabla_{
x}\tau\>\] for any vector field $v$ and all 1-forms $\tau$.

\begin{propos}\label{tor} For a compatible preconnection in
the sense of (\ref{eqq1})-(\ref{eqq2}), the torsion   obeys
\[ T(\hat x,\hat y)(\extd z)+{\rm cyclic}=0,\quad \forall x,y,z\in
C^\infty(M).\]
\end{propos}
\proof From the definitions, we have
\begin{eqnarray}\label{tor1}
  \<T(\hat x,\hat
y),\extd z\>=\<\hat x,\hat\nabla_{ y}\extd z\>-\<\hat
y,\hat\nabla_{ x}\extd z\>
\end{eqnarray}
  for all $x,y,z\in C^\infty(M)$, using the definition
of torsion, converting to an operation on forms and using the
Jacobi identity for the Poisson bracket. We then take the cyclic
sum $x\to y\to z\to x$, use (\ref{eqq2}) three times, and the
Jacobi identity again.  \eproof

This completes our general comments. For further explicit
computations, we will suppose that $M$ is finite-dimensional with
coordinate patch with coordinate functions $(x^1,\dots,x^n)$. We
use the summation convention for repeated indices. Suppose that
the Poisson structure is given by
\[
\{y,z\}\,=\,\omega^{ij}\,\frac{\del y}{\del x^i}
\,\frac{\del z}{\del x^j}\ ,
\]
where $\omega^{ij}\,=\,-\,\omega^{ji}$ and $y,z$ are functions.
The Jacobi identity for $\{\ ,\ \}$ gives
\begin{eqnarray}\label{projacobi}
\omega^{is}\,\frac{\del \omega^{jk}}{\del x^s}\,+\,
\omega^{js}\,\frac{\del \omega^{ki}}{\del x^s}\,+\,
\omega^{ks}\,\frac{\del \omega^{ij}}{\del x^s}\,=\,0\ .
\end{eqnarray} From (\ref{eqq3}) we can write
\[
\gamma(y,\extd x^i)\,=\, c_n^{\phantom{j}ik}\, \frac{\del
y}{\del x^k}\, \extd x^n\ ,
\]
or alternatively $\gamma(x^k,\extd
x^i)\,=\,c_n^{\phantom{j}ik}\,\extd x^n$. Then from (\ref{eqq1})
we get
\[
\gamma(y,a_i\,\extd x^i)\,=\, \Big(\omega^{kq}\, \frac{\del
a_n}{\del x^q}\,+\, c_n^{\phantom{j}ik}\,a_i\Big)\,
\frac{\del y}{\del x^k}\, \extd x^n\ .
\]
From (\ref{eqq2}) we get
\begin{eqnarray}\label{proomega}
\frac{\del \omega^{ij}}{\del
x^n}\,=\,c_n^{\phantom{n}ji}\,-\, c_n^{\phantom{n}ij}\ .
\end{eqnarray}
Given a function $y$ the associated Hamiltonian vector field $\hat
y$ is of course defined by
\[
\hat y\,=\,\omega^{kq}\, \frac{\del y}{\del x^k}\,
  \frac{\del }{\del x^q}\ .
\]

\subsection{The symplectic case}

To simplify the computations we now specialise to the
non-degenerate case,  where the matrix $\omega^{ij}$ is
invertible. We write its inverse as $\omega_{ij}$, so that
$\omega^{ij}\,\omega_{jk}\,=\,\delta_k^i$. Then equation
(\ref{projacobi}) shows that the 2-form $\omega_{ij}\,\extd x^i
\wedge \extd x^j$ is closed, and so is a symplectic form. We now
use material from \cite{GRS-Fed}. We can rewrite the formula for
$\gamma$ as
\[
\gamma(y,\tau)\,=\, \omega^{kq}\, \frac{\del y}{\del
x^k}\, \nabla_q(\tau)\ ,
\]
where $\nabla$ is now a fully defined covariant connection,
\[
\nabla_q(a_n\, \extd x^n)\,=\,\Big(\frac{\del a_n}{\del
x^q}\,+\, \omega_{qs}\, c_n^{\phantom{n}is}\, a_i\Big)\, \extd
x^n\,=\, \Big(\frac{\del a_n}{\del x^q}\,-\,
\Gamma_{qn}^i\, a_i\Big)\, \extd x^n\ .
\]
The formula for the Christoffel symbols is
\[
\Gamma_{qn}^i\,=\, -\,\omega_{qs}\, c_n^{\phantom{n}is} \quad{\rm
or}\quad c_n^{\phantom{n}ik}\,=\,-\,\omega^{kq}\,\Gamma_{qn}^i \ .
\]

From (\ref{proomega}) we see that
\begin{eqnarray}\label{proomega2}
\frac{\del \omega^{ij}}{\del x^n}\,+\,
\omega^{iq}\,\Gamma_{qn}^j\,+\, \omega^{qj}\,\Gamma_{qn}^i\,=\,0\
, \end{eqnarray} which can be rewritten using the torsion tensor
$T_{qn}^j=\Gamma_{qn}^j-\Gamma_{nq}^j$ as
\begin{eqnarray}\label{proomega3}
\nabla_n \omega^{ij}\,+\, \omega^{iq}\,T_{qn}^j\,+\,
\omega^{qj}\,T_{qn}^i\,=\,0\ . \end{eqnarray}

\begin{propos}\label{omegacen} The 2-form $\omega$ commutes in the
$\bullet$ product with all functions to $O(\hbar^{2})$ if and only
if the connection preserves $\omega$ and is torsion free.
\end{propos} \proof We use Darboux coordinates where $\omega_{ij}$
is constant. If the 2-form $\omega$ is central, then we see that
the connection preserves $\omega$, i.e.\
\begin{eqnarray*}
\omega^{iq}\,\Gamma_{nq}^j\,+\, \omega^{qj}\,\Gamma_{nq}^i\,=\,0\
, \end{eqnarray*} From (\ref{proomega2}) we see that
\begin{eqnarray*}
\omega^{iq}\,\Gamma_{qn}^j\,+\, \omega^{qj}\,\Gamma_{qn}^i\,=\,0\
, \end{eqnarray*} As in \cite{GRS-Fed} we set
$\Gamma_{ikj}=\omega_{il}\,\Gamma^{l}_{kj}$, and then the last two
equations give $\Gamma_{snr}\,=\,\Gamma_{rns}$ and
$\Gamma_{srn}\,=\,\Gamma_{rsn}$, so $\Gamma_{snr}$ is totally
symmetric and the torsion vanishes. \eproof

\subsection{The second Jacobi identity in the symplectic case}

\begin{lemma}
     For a 1-form $\tau$ and a vector field $Y$,
\begin{eqnarray*}
     \extd (\nabla_Y \tau)-\nabla_Y \extd \tau &=&(\,
(T^n_{jl}\,Y^{j}+Y^{n}_{\phantom{n};l})\,\tau_{i;n}-N^k_{jil}\,
\tau_k\,Y^{j} \,)\extd x^l\wedge \extd x^i\ , \end{eqnarray*}
where
\[
N^k_{jil}\,=\, T^k_{ji;l}\,+\,
\frac12\,R^{k}_{\phantom{k}jli}\,+\, T^n_{jl}\,T^k_{in}\,+
\,\frac12\,T^n_{li}\,T^k_{jn}\ .
\]
\end{lemma}
   \proof
   \begin{eqnarray*}
\extd (\nabla_Y \tau)-\nabla_Y \extd \tau &=&(\,
(\Gamma^n_{jl}\,\tau_{i,n}-\Gamma^k_{ji,l}\,\tau_k)Y^j\,+\,
(\tau_{i,j}-\Gamma^k_{ji}\,\tau_k)Y^j_{,l}\,)\extd x^l\wedge \extd
x^i \cr &=&(\, (\Gamma^n_{jl}\,\tau_{i,n}-T^k_{ji,l}\,\tau_k
-\Gamma^k_{ij,l}\,\tau_k)Y^j\,+\, \tau_{i;j}\,Y^j_{,l}\,)\extd
x^l\wedge \extd x^i \end{eqnarray*} Now we use
\begin{eqnarray*}
R^{k}_{\phantom{k}jli}\,\extd x^l\wedge \extd x^i &=&(
\Gamma^k_{ij,l}\,-\, \Gamma^k_{lj,i}\,+\,
\Gamma^m_{ij}\,\Gamma^k_{lm}\,-\,\Gamma^m_{lj}\,\Gamma^k_{im}) \,
\extd x^l\wedge \extd x^i \cr &=& 2\,(\Gamma^k_{ij,l}\,+\,
\Gamma^m_{ij}\,\Gamma^k_{lm}) \, \extd x^l\wedge \extd x^i
\end{eqnarray*} to get
\begin{eqnarray*}
\extd (\nabla_Y \tau)-\nabla_Y \extd \tau &=&(\,
(\Gamma^n_{jl}\,\tau_{i,n}-T^k_{ji,l}\,\tau_k
-\frac12\,R^{k}_{\phantom{k}jli}\,\,\tau_k +
\Gamma^n_{ij}\,\Gamma^k_{ln}\,\tau_k)Y^j\,+\,
\tau_{i;j}\,Y^j_{,l}\,)\extd x^l\wedge \extd x^i \cr &=&(\,
(\Gamma^n_{jl}\,\tau_{i,n}-T^k_{ji,l}\,\tau_k
-\frac12\,R^{k}_{\phantom{k}jli}\,\,\tau_k +
\Gamma^n_{ij}\,\Gamma^k_{ln}\,\tau_k-\Gamma^n_{lj}\,(\nabla_n\tau_i))Y^j
\cr &&\,+\,(\nabla_j\tau_i)\,(\nabla_l Y^j) \,)\extd x^l\wedge
\extd x^i \cr &=&(\, (T^n_{jl}\,\tau_{i;n}-T^k_{ji,l}\,\tau_k
-\frac12\,R^{k}_{\phantom{k}jli}\,\,\tau_k +
\Gamma^n_{ij}\,\Gamma^k_{ln}\,\tau_k+\Gamma^n_{jl}\,\Gamma^k_{ni}\,
\tau_k)Y^j \cr &&\,+\,(\nabla_j\tau_i)\,(\nabla_l Y^j) \,)\extd
x^l\wedge \extd x^i \cr &=&(\,
(T^n_{jl}\,\tau_{i;n}-T^k_{ji;l}\,\tau_k
-\frac12\,R^{k}_{\phantom{k}jli}\,\,\tau_k +
\Gamma^n_{ji}\,\Gamma^k_{ln}\,\tau_k+\Gamma^n_{jl}\,\Gamma^k_{ni}\,\tau_k
\cr &&
-\,\Gamma^n_{lj}\,T^k_{ni}\,\tau_k-\Gamma^n_{li}\,T^k_{jn}\,\tau_k
)Y^j \,+\,(\nabla_j\tau_i)\,(\nabla_l Y^j) \,)\extd x^l\wedge
\extd x^i \cr &=&(\, (T^n_{jl}\,\tau_{i;n}-T^k_{ji;l}\,\tau_k
-\frac12\,R^{k}_{\phantom{k}jli}\,\,\tau_k -
T^n_{jl}\,T^k_{in}\,\tau_k-\frac12\,T^n_{li}\,T^k_{jn} \,\tau_k
)Y^j \cr &&\,+\,(\nabla_j\tau_i)\,(\nabla_l Y^j) \,)\extd
x^l\wedge \extd x^i \cr &=&(\,
(T^n_{jl}\,Y^{j}+Y^{n}_{\phantom{j};l})\,\tau_{i;n}-N^k_{jil}\,\tau_k\,Y^{j}
\,)\extd x^l\wedge \extd x^i\ .\square \end{eqnarray*}

\begin{propos} The super commutator between two 1-forms $\tau$ and
$\xi$ is
\[
[\xi,\tau]_{\bullet}\,=\, \hbar\,\omega^{jn}\,\nabla_{j}\xi\wedge
\nabla_{n}\tau\,+\,\hbar\,\omega^{jn}\,N^{k}_{jil}\,\xi_{n}\,\tau_{k}
\,\extd x^l\wedge \extd x^i\,+\,O(\hbar^{2})\ .
\]
By the symmetry of $[\xi,\tau]_{\bullet}$ we must have
$\omega^{jn}\,N^{k}_{jil}-\omega^{jn}\,N^{k}_{jli}$ symmetric in
$nk$.
\end{propos} \proof It is enough to consider $\xi=a\,\extd b$.  Begin
with
\begin{eqnarray*}
     [\extd b,\tau]_{\bullet} &=&
\extd \,[b,\tau]_{\bullet}\,-\,[b,\extd \tau]_{\bullet}
\,+\,O(\hbar^{2}) \cr &=& \hbar\,(\extd \,\nabla_{\hat
b}-\nabla_{\hat b}\,\extd )\tau \,+\,O(\hbar^{2}) \cr &=& \hbar\,
(\, (T^n_{jl}\,{\hat b}^{j}+{\hat b}^{n}_{\phantom{n};l})
\,\tau_{i;n}-N^k_{jil}\,\tau_k\,{\hat b}^{j} \,)\extd x^l\wedge
\extd x^i \,+\,O(\hbar^{2})\ .
\end{eqnarray*} Now we use $\hat b^{n}=\omega^{in}\,b_{,i}$ to
get
\begin{eqnarray*}
     [\extd b,\tau]_{\bullet} &=&
     \hbar\,
(\, (T^n_{jl}\,\omega^{mj}\,b_{,m}+
\omega^{mn}_{\phantom{mn};l}\,b_{,m}+\omega^{mn}\,(b_{,m})_{;l} )
\,\tau_{i;n}\\
&&-N^k_{jil}\,\tau_k\,\omega^{mj}\,b_{,m} \,)\extd x^l\wedge \extd
x^i \,+\,O(\hbar^{2})\ . \end{eqnarray*} Using equation
(\ref{proomega3}) for the covariant derivative of $\omega$,
\begin{eqnarray*}
     [\extd b,\tau]_{\bullet} &=&
     \hbar\,
(\, (T^m_{jl}\,\omega^{nj}\,b_{,m}+\omega^{mn}\,(b_{,m})_{;l} )
\,\tau_{i;n}-N^k_{jil}\,\tau_k\,\omega^{mj}\,b_{,m} \,)\extd
x^l\wedge \extd x^i \,+\,O(\hbar^{2})\cr &=&
     \hbar\,
(\, (T^m_{jl}\,\omega^{nj}\,b_{,m}+\omega^{jn}\,(b_{,j})_{;l} )
\,\tau_{i;n}-N^k_{jil}\,\tau_k\,\omega^{nj}\,b_{,n} \,)\extd
x^l\wedge \extd x^i \,+\,O(\hbar^{2})\cr &=&
     \hbar\,\omega^{nj}\,
(\, (T^m_{jl}\,b_{,m}-\,(b_{,j})_{;l} )
\,\tau_{i;n}-N^k_{jil}\,\tau_k\,b_{,n} \,)\extd x^l\wedge \extd
x^i \,+\,O(\hbar^{2})\cr &=&
     \hbar\,\omega^{nj}\,
(\, -\,(b_{,l})_{;j} \, \,\tau_{i;n}-N^k_{jil}\,\tau_k\,b_{,n}
\,)\extd x^l\wedge \extd x^i \,+\,O(\hbar^{2}) \ . \end{eqnarray*}
Now using the next equation gives the answer;
\begin{eqnarray*}
     [a\,\extd b,\tau]_{\bullet} &=& a\,[\extd b,\tau]_{\bullet} \,+\,
\extd b\wedge [a,\tau]_{\bullet}\,+\,O(\hbar^{2}) \ .\quad\square
\end{eqnarray*}

\begin{propos} If the curvature vanishes and
$E^{nk}_{\phantom{nk}li}\equiv \omega^{jn}\,N^{k}_{jil}$ is
covariantly constant, then the obstruction $J_2$ to the second
super Jacobi identity vanishes.  \end{propos} \proof First
consider the case $E=0$.  In Darboux coordinates where $\omega$ is
constant,
\begin{eqnarray*}
  &&   \nquad [a,[\xi,\tau]_{\bullet}]_{\bullet}\,-\,
[\tau,[a,\xi]_{\bullet}]_{\bullet}\,+\,
[\xi,[\tau,a]_{\bullet}]_{\bullet} \cr &=& \hbar(\,
[a,\omega^{jn}\,\nabla_{j}\xi\wedge\nabla_{n}\tau]_{\bullet}\,-\,
[\tau,\nabla_{\hat a}\xi]_{\bullet}\,-\, [\xi,\nabla_{\hat
a}\tau]_{\bullet}\,)\,+\,O(\hbar^{3}) \cr &=& \hbar^{2}\,(\,
\omega^{jn}\,[\nabla_{\hat a},\nabla_{j}]\xi\wedge\nabla_{n}\tau
\,+\, \omega^{jn}\,\nabla_{j}\xi\wedge[\nabla_{\hat
a},\nabla_{n}]\tau \,)\,+\,O(\hbar^{3}) \ . \end{eqnarray*} If the
curvature is zero this becomes
\begin{eqnarray*}
&&\nquad \hbar^{2}\,(\, \omega^{jn}\,\nabla_{\{\hat
a,j\}}\xi\wedge\nabla_{n}\tau \,+\,
\omega^{jn}\,\nabla_{j}\xi\wedge\nabla_{\{\hat a,n\}}\tau
\,)\,+\,O(\hbar^{3}) \cr &=& \hbar^{2}\,(\,
\omega^{mn}\,\nabla_{\{\hat a,m\}}\xi\wedge\nabla_{n}\tau \,+\,
\omega^{jm}\,\nabla_{j}\xi\wedge\nabla_{\{\hat a,m\}}\tau
\,)\,+\,O(\hbar^{3}) \cr &=& \hbar^{2}\,(\, \omega^{mn}\,\{\hat
a,\del_{m}\}^{j} \,+\, \omega^{jm}\,\{\hat
a,\del_{m}\}^{n}
\,)\,\nabla_{j}\xi\wedge\nabla_{n}\tau\,+\,O(\hbar^{3})\ .
\end{eqnarray*} Now we calculate
\begin{eqnarray*} \omega^{mn}\,\{\hat a,\del_{m}\}^{j}
     &=& -\,\omega^{mn}\,(\hat a^{j})_{,m}\,=\,
     \omega^{nm}\,\omega^{ij}\,a_{,im}\,=\,
      \omega^{nm}\,\omega^{ij}\,a_{,mi}\cr
      &=& \omega^{ni}\,\omega^{mj}\,a_{,im}\,=\,
      -\,\omega^{jm}\,\omega^{ni}\,a_{,im}\,=\,-\,
      \omega^{jm}\,\{\hat a,\del_{m}\}^{n}
      \ .
\end{eqnarray*} Now we consider the general case, where
\begin{eqnarray*}
  &&\nquad [a,[\xi,\tau]_{\bullet}]_{\bullet}\,-\,
[\tau,[a,\xi]_{\bullet}]_{\bullet}\,+\,
[\xi,[\tau,a]_{\bullet}]_{\bullet}\\ &=& \hbar\,
[a,E^{nk}_{\phantom{nk}li}\xi_{n}\tau_{k}\,\extd x^l\wedge \extd
x^i]_{\bullet}\,-\, E^{nk}_{\phantom{nk}li}\tau_{n}(\nabla_{\hat
a}\xi)_{k}\, \extd x^l\wedge \extd x^i\,\\
&&-\, E^{nk}_{\phantom{nk}li}\xi_{n}(\nabla_{\hat a}\tau)_{k}\,
\extd x^l\wedge \extd x^i \,+\,O(\hbar^{3})\ , \end{eqnarray*}
which leads to $\nabla_{{\hat a}}E^{nk}_{\phantom{nk}li}=0$.
 \eproof

\subsection{The third Jacobi identity in the symplectic case}

We already know from Section~2 that if $J_1=J_2=0$ then the third
part $J_3$ of the super Jacobi identity, for three 1-forms,
vanishes. In view of the above, this is the case when the
curvature and torsion both vanish.  Here we give, as a check, a
direct proof of this fact. We consider the super Jacobi identity
with three 1-forms, $\phi$, $\tau$ and $\eta$:
\begin{eqnarray*}
[\phi,[\tau,\eta]_\bullet]_\bullet &=& [\phi,
\hbar\,\omega^{kq}\,\nabla_k\tau\wedge\nabla_q\eta
\,+\,O(\hbar^2)]_\bullet \cr
  &=& \hbar\,[\phi,
\omega^{kq}\,\nabla_k\tau\wedge\nabla_q\eta ]_\bullet
\,+\,O(\hbar^3) \cr &=& \hbar\,\omega^{kq}\, ([\phi,
\nabla_k\tau]_\bullet\wedge\nabla_q\eta -
\nabla_k\tau\wedge[\phi,\nabla_q\eta]_\bullet)
  \,+\,O(\hbar^3)\cr
&=& \hbar^2\,\omega^{kq}\,\omega^{nm}\, (\nabla_n\phi\wedge
\nabla_m\nabla_k\tau\wedge\nabla_q\eta -
\nabla_k\tau\wedge\nabla_n\phi\wedge\nabla_m\nabla_q\eta)
  \,+\,O(\hbar^3)\cr
&=& \hbar^2\,\omega^{kq}\,\omega^{nm}\, (\nabla_n\phi\wedge
\nabla_m\nabla_k\tau\wedge\nabla_q\eta +
\nabla_k\tau\wedge\nabla_m\nabla_q\eta\wedge\nabla_n\phi)
  \,+\,O(\hbar^3)
\ . \end{eqnarray*} To show that the cyclic sum vanishes to
$O(\hbar^2)$ we need the following 3-form to be zero:
\begin{eqnarray*}
&&\nquad \omega^{kq}\,\omega^{nm}\, ( \nabla_n\phi\wedge
\nabla_m\nabla_k\tau\wedge\nabla_q\eta +
\nabla_k\tau\wedge\nabla_m\nabla_q\eta\wedge\nabla_n\phi \cr &&
\nquad + \nabla_n\tau\wedge \nabla_m\nabla_k\eta\wedge\nabla_q\phi
+ \nabla_k\eta\wedge\nabla_m\nabla_q\phi\wedge\nabla_n\tau  \cr &&
\nquad + \nabla_n\eta\wedge \nabla_m\nabla_k\phi\wedge\nabla_q\tau
+ \nabla_k\phi\wedge\nabla_m\nabla_q\tau\wedge\nabla_n\eta )   \cr
&=& (\omega^{kq}\,\omega^{nm}+\omega^{nk}\,\omega^{qm})\,
\nabla_n\phi\wedge \nabla_m\nabla_k\tau\wedge\nabla_q\eta   \cr &&
+\, (\omega^{kq}\,\omega^{nm}+\omega^{qn}\,\omega^{km})\,
\nabla_k\tau\wedge\nabla_m\nabla_q\eta\wedge\nabla_n\phi\cr && +\,
(\omega^{kq}\,\omega^{nm}+\omega^{qn}\,\omega^{km})\,
\nabla_k\eta\wedge\nabla_m\nabla_q\phi\wedge\nabla_n\tau\cr &=&
(\omega^{kq}\,\omega^{nm}+\omega^{nm}\,\omega^{qk})\,
\nabla_n\phi\wedge \nabla_m\nabla_k\tau\wedge\nabla_q\eta   \cr &&
+\, (\omega^{kq}\,\omega^{nm}+\omega^{mn}\,\omega^{kq})\,
\nabla_k\tau\wedge\nabla_m\nabla_q\eta\wedge\nabla_n\phi\cr && +\,
(\omega^{kq}\,\omega^{nm}+\omega^{mn}\,\omega^{kq})\,
\nabla_k\eta\wedge\nabla_m\nabla_q\phi\wedge\nabla_n\tau\cr &=&0\
.\quad\square \end{eqnarray*}

\subsection{The generalised braiding associated to $\nabla$}

If $A$ is an algebra equipped with a calculus and a left covariant
derivative $\nabla:\Omega^1\to \Omega^1\tens_A\Omega^1$ on it
(meaning $\nabla(a\bullet\tau)=\extd a\tens_A\tau+a\bullet
\nabla(\tau)$ for all $a\in A,\tau\in \Omega^1$) it may be
possible to define a `generalised braiding'\cite{Madore}:
\[ \bar\sigma:\Omega^1\tens_A\Omega^1\to \Omega^1\tens_A
\Omega^1,\quad \sigma(\tau\tens \extd a)=\nabla(\tau\bullet
a)-\nabla(\tau)\bullet a,\] where $\sigma$ is given on exact forms
in the second argument and assumes that $\ker\extd$ is spanned by
$1$. It is already a left module map and in well-behaved cases it
extends also as a right module map to $\Omega^1\tens \Omega^1$, in
which case it follows using the Leibniz rule that it descends to
an operator $\bar\sigma$ on $\Omega^1\tens_A \Omega^1$.

In our case $\bullet$ is the deformed module structure on
$\Omega^1(M)$, etc. and we suppose a connection coinciding at
lowest order with a given connection (such as the symplectic one
above). Then \eqn{sigma}{ \sigma(\tau\tens \extd x)\,=\, \extd
x\bar\tens\tau\,-\,\nabla([x,\tau]_\bullet)\,+\,[x,\nabla\tau]_\bullet\
=\extd x\bar\tens\tau+\hbar\rho(\extd x\tens \tau)+O(\hbar^2)}
where we further suppose an expansion as shown. We write
$\bar\tens$ for the tensor product over $C^\infty(M)$.

\begin{propos} $\sigma$ defining a generalised braiding corresponds at
the semiclassical level to $\rho$ obeying
\[
\eta\bar\tens\gamma(x,\tau)+\rho(x\eta\tens\tau)
=x\rho(\eta\tens\tau)=\rho(\eta\tens \tau
x)-\gamma(x,\eta)\bar\tens\tau\] for all $x\in C^\infty(M)$ and
$\tau,\eta\in \Omega^1(M)$. Moreover, when $\nabla$ is given by
our symplectic connection in Section~3.1, we have
\begin{eqnarray*}
\rho(\eta\tens\tau) &=&\omega^{jq}\,\omega^{is} \,\eta_j\,\tau_i\,
R_{smkq} \,\extd x^k\bar\tens \extd
x^m\,-\,\omega^{sq}\,\nabla_s\eta\bar \tens \nabla_q \tau
\end{eqnarray*}
in a coordinate chart where $\nabla=\extd x^k\bar\tens\nabla_k$.
\end{propos}
\proof We suppose similarly that
$\sigma(\eta\tens\tau)=\tau\bar\tens\eta+\hbar\rho(\tau\tens\eta)+O(\hbar)$.
That $\sigma$ is a left module map translates to $\tau\bar\tens
a\bullet\eta+\hbar\rho(\tau\tens
a\bullet\eta)+O(\hbar^2)=a\bullet\tau\bar\tens\eta+\hbar
a\bullet\rho(\tau\tens\eta)+O(\hbar^2)$. Since
$[a,\tau]_\bullet=\hbar\gamma(a,\tau) + O(\hbar^2)$ and $a$ passes
through $\bar\tens$, and since everything commutes at lowest
order, we obtain the first equality. Similarly, that $\sigma$ is a
right module map implies the second equality. That $\sigma$
descends to $\bar\tens$ gives nothing new, appearing as equality
of the first and last expressions.

Moreover, we compute
\begin{eqnarray*}
\rho( \extd x^j\tens\tau ) &=&  -\extd x^k\bar\tens
\nabla_k\gamma(x^j,\tau)\,+\,\gamma(x^j, \extd x^k\bar\tens
\nabla_k(\tau )) \cr &=& -\,\extd x^k\bar\tens
(\nabla_k\gamma(x^j,\tau )-\gamma(x^j,\nabla_k(\tau)))+\,
\gamma(x^j, \extd x^k)\bar \tens \nabla_k(\tau )\ .
\end{eqnarray*} Now set $X=\omega^{jq}\del_q$, and then
\begin{eqnarray*}
\rho(\extd x^j\tens\tau) &=& -\extd x^k\bar\tens
(\nabla_k\nabla_X-\nabla_X\nabla_k)(\tau ) \,+\, \nabla_X( \extd
x^k)\bar\tens \nabla_k(\tau )\ .
\end{eqnarray*}
In Darboux coordinates $[\del_k,X]=0$, and
\begin{eqnarray*}
\rho(\extd x^j\tens\tau) &=& -\extd x^k\bar\tens
R(\del_k,X)(\extd x^i)\, +\,\nabla_X( \extd x^k)\bar\tens
\nabla_k(\tau )\ .
\end{eqnarray*} Now we calculate the terms as
\begin{eqnarray*}
\extd x^k\bar\tens R(\del_k,X)(\tau ) &=& \omega^{jq}\, \extd
x^k\bar\tens R(\del_k,\del_q)(\tau ) \,=\,
-\,\omega^{jq}\,\tau_i\, \extd x^k \bar\tens
R^i_{\phantom{i}mkq}\, \extd x^m \ ,\cr \nabla_X( \extd
x^k)\bar\tens \nabla_k(\tau ) &=& \omega^{jq}\nabla_q( \extd
x^k)\bar\tens \nabla_k(\tau ) \ .
\end{eqnarray*}
Using these, from (\ref{proomega2}) we have
\begin{eqnarray*}
\rho(\extd x^j\tens\tau) &=&\omega^{jq}\,\tau_i\,
R^i_{\phantom{i}mkq} \,\extd x^k\bar\tens \extd
x^m-\,\omega^{sq}\,\nabla_s(\extd x^j)\bar \tens \nabla_q \tau.
\end{eqnarray*}
Finally, we deduce the general form of $\rho$ from its properties
above as
\begin{eqnarray*}
\rho(\eta_j\extd x^j \tens\tau) &=&
   \omega^{jq}\,\tau_i\, \eta_j\,
R^i_{\phantom{i}mkq} \,\extd x^k\bar\tens \extd
x^m\,-\,\omega^{sq}\,\eta_j\,\nabla_s(\extd x^j) \bar\tens
\nabla_q \tau\,-\,   \extd x^j\tens\gamma(\eta_j,\tau).
\end{eqnarray*}
Putting $\eta=\eta_j\extd x^j$ and assuming that $\nabla$ is the
symplectic connection defined by $\gamma$ in Section~3.1 gives the
result stated. We use the notations from \cite{GRS-Fed}. \eproof

It is also instructive to look at the braid relations, which
justifies the term `generalised braiding'. Writing
$\bar\sigma_{12}$ for $\bar\sigma$ acting in the first and second
tensor factors, etc, the expression
\begin{eqnarray*}
(\bar\sigma_{12}\bar\sigma_{23}\bar\sigma_{12}
-\bar\sigma_{23}\bar\sigma_{12}\bar\sigma_{23})(\xi\bar\tens\eta\bar\tens\tau)
\end{eqnarray*}
can be calculated to leading order in $\hbar$ as the quotient in
$(\Omega^1)^{\bar\tens 3}$ of
\begin{eqnarray*}
\big(\tilde\rho_{12}\,\tilde\rho_{23}+\tilde\rho_{13}\,\tilde\rho_{23}+\tilde\rho_{12}\,\tilde\rho_{13}
-\tilde\rho_{13}\,\tilde\rho_{12}-\tilde\rho_{23}\,\tilde\rho_{12}-\tilde\rho_{23}\,\tilde\rho_{13}\big)
(\tau\tens\eta\tens\xi)\ ,
\end{eqnarray*}
where $\tilde\rho:\Omega^1\tens \Omega^1\to \Omega^1\tens\Omega^1$
is any lift of the output of $\rho$. The left hand side is the
`classical Yang-Baxter expression'. In our case the natural lift
determined by our coordinate system of the $\rho$ in
Proposition~3.7 is
\begin{eqnarray}\label{tilderho}
&&\tilde\rho(\tau\tens\eta) \,= \omega^{jq}\,\omega^{is}
\,\tau_j\,\eta_i\, R_{smkq} \,\extd x^k\tens \extd
x^m\,-\,\omega^{sq}\,\nabla_s\tau \tens \nabla_q \eta.
\end{eqnarray}
On substituting this into the classical Yang-Baxter expression,
the part coming from the $\omega^{sq}\,\nabla_s\eta \tens \nabla_q
\xi$ term is
\begin{eqnarray*}
&& \omega^{sq}\,\omega^{ij}\,\big(\nabla_s\tens[\nabla_q,\nabla_i]
\tens\nabla_j +\nabla_s\tens\nabla_i\tens[\nabla_q,\nabla_j]
+[\nabla_s,\nabla_i]\tens\nabla_q\tens\nabla_j\big)
\end{eqnarray*}
applied to $\zeta\tens\eta\tens\xi$. We deduce that if the
curvature vanishes, then the classical Yang-Baxter equation holds
for $\tilde\rho$ (even before we quotient), and hence $\bar\sigma$
satisfies the braid relation to $O(\hbar^3)$. We can also
  calculate
\begin{eqnarray*}
(\tilde\rho+\tilde\rho_{21})(\tau\tens\eta) &=&  \omega^{jq}\,
\omega^{is}\,R_{smkq} (\tau_j\,\eta_i\, \extd x^k\tens\extd
x^m+\eta_j\,\tau_i\, \extd x^m\tens\extd x^k)\cr &&\,
-\,\omega^{sq}(\nabla_s\tau\tens\nabla_q\eta+\nabla_q\tau
\tens\nabla_s\eta) \cr &=&
  \omega^{jq}\,\omega^{is}\,(R_{smkq}+R_{qkms})\,\tau_j\,\eta_i\,
\extd x^k\tens\extd x^m
\end{eqnarray*}
so if the curvature vanishes, then $\tilde\rho$ is antisymmetric,
and $\bar\sigma^2$ is the identity to $O(\hbar^2)$.

\subsection{Example: the noncommutative torus}

Here we give an elementary example of the above, namely the
noncommutative torus in an algebraic form generated by invertible
$u,v$ with $vu=quv$. The natural noncommutative calculus here is
generated by
\[ \tau^1=u^{-1}\extd u,\quad \tau^2=v^{-1}\extd v\]
with relations
\[ [x,\tau^i]=0,\quad \tau^1\wedge\tau^2+\tau^2\wedge\tau^1=0,\quad
\tau^i\wedge\tau^i=0\] for all functions $x$ (this is not as
trivial as it seems, it implies $v\extd u=q(\extd u)v$, etc.)

For our purposes we take $q=e^{-\hbar}$ as a formal expansion
parameter, which we stress is only the `trivial' part of the
theory from the usual point of view (where one would have
$\hbar=2\pi\imath\theta$ with $\theta\in[0,1]$ and complete to a
$C^*$-algebra). Next, in the classical limit we identify
$u=e^{\imath\theta^1}$ and $v=e^{\imath\theta^2}$ in terms of
local angle coordinates on $S^1\times S^1$. Along with their
inverses they generate a subalgebra inside $C^\infty(S^1\times
S^2)$ and induce a formal deformation of the whole algebra via a
Moyal-type product. For our purposes all structures can be
extracted just from the subalgebra with the above generators. Note
that $\tau^i=\imath\extd\theta^i$.

We start with the Poisson structure which is clearly
\[ \{u,v\}=1=-\{v,u\},\quad \{u,u\}=\{v,v\}=0\] which means
$(\omega^{ij})=\begin{pmatrix}0&1\\-1&0\end{pmatrix}$ in the
$\tau^i$ basis. To see this, note that
\[ \extd x=(\del_i x)\tau^i;\quad \del_1=u{\del\over\del u},\quad
\del_2=v{\del\over\del v}\ .\] Then $\{u,y\}=\del_2 y$ and
$\{v,y\}=-\del_1 y$ for all $y$ since $\hat u$ and $\hat v$ extend
as derivations. Then since each $\hat y$ is a derivation, we
similarly deduce
\[ \{x,y\}=\omega^{ij}(\del_i x)(\del_j y)\]
with $\omega^{ij}$ as above. The symplectic form is
$\tau^1\wedge\tau^2$ up to normalisation.

Finally, since $\tau^i$ are central, we know that
$\gamma(x,\tau^i)=0$ for  $i=1,2$ and all functions $x$. Hence
\begin{eqnarray} \label{haaa}
\nabla_{\hat
x}(\eta_i\tau^i)=\gamma(x,\eta_i\tau^i)=\{x,\eta_i\}\tau^i
\end{eqnarray}
by the connection property. This is our symplectic connection.
Compatibility in the sense of (\ref{eqq2}) holds as it must,
reducing to \[ \del_i\{x,y\}=\{\del_ix,y\}+\{x,\del_i y\}\] which
holds because $\omega^{ij}$ are constant. We may also compute its
curvature:
\begin{eqnarray*}
R(\hat x,\hat y)a_i\tau^i&=&\nabla_{\hat
x}\{y,a_i\}\tau^i-\nabla_{\hat
y}\{x,a_i\}\tau^i-\nabla_{\widehat{\{x,y\}}}a_i\tau^i\\
&=&\{x,\{y,a_i\}\}\tau^i-\{y,\{x,a_i\}\}\tau^i-\{\{x,y\},a_i\}\tau^i=0
\end{eqnarray*}
as it must since the above quantum differential calculus is
associative to all orders. Moreover,
\begin{eqnarray*} \<T(\hat x,\hat y),\extd u\>
&=&\<\hat x,\nabla_{\hat y}\extd
u\>-\<\hat y,\nabla_{\hat x}\extd u\>\\
&=&\<\hat x,\{y,\del_i u\}\tau^i\> -\<\hat y,\{x,\del_i
u\}\tau^i\>\\
&=&\{y,u\}\<\hat x,\tau^1\>-\{x,u\}\<\hat y,\tau^1\>\\
&=&\{y,u\}u^{-1}\{x,u\}-\{x,u\}u^{-1}\{y,u\}=0
\end{eqnarray*}
since $\<\hat x,u^{-1}\extd u\>=u^{-1}\hat x(u)=u^{-1}\{x,u\}$.
Similarly $\<T(\hat x,\hat y),\extd v\>=0$, and hence $T(\hat
x,\hat y)=0$. Hence the torsion vanishes as it must by
Proposition~\ref{omegacen}. So we have a flat torsion free (and
symplectic) connection on the standard torus induced by
noncommutative geometry.

Finally, the semiclassical braiding in Proposition~3.7 for this
connection, using the fact that the curvature is zero, and
(\ref{haaa}), comes out to be $\rho(\tau^i\tens \tau^j)=0$ so that
\[ \rho(\eta_i\tau^i\tens \xi_j\tau^j)=-\{\eta_i,\xi_j\}\tau^i\bar\tens\tau^j\]
on general 1-forms with coefficients $\eta_i,\xi_j$ in the
left-invariant basis. Note also that
$\bar\sigma(\tau_i\bar\tens\tau_j)=
\tau_j\bar\tens\tau_i+O(\hbar^2)$ agrees with the quantum
case\cite{EJB-Vec} where it was found that the unique braiding on
the noncommutative torus which is compatible with interior
products for the given differential structure is the usual flip on
the $\tau^i$.

\section{Application to Poisson-Lie groups}

In this section we are going to analyse the meaning of the
deformation results of Section~2 in a different geometrical
setting, namely that of Poisson-Lie groups. This is at the other
extreme from the previous section in that these are never
symplectic (the Poisson bracket being degenerate at the group
identity). On the other hand one can impose the group symmetry and
reduce expressions to the Lie algebra level which has its own
interest.

Poisson-Lie groups are the semiclassical data for Hopf algebras or
quantum groups. At the Lie level they correspond to a Lie
bialgebra $(\cg,\delta)$ where $\delta:\cg\to \cg\tens \cg$ is a
Lie cobracket making $\cg^*$ a Lie algebra and forming a cocycle
in $Z^1_\ad(\cg,\cg\tens \cg)$. This is like an infinitesimal
quantum group. Exponentiating the Lie algebra to a group $G$ and
$\delta$ to $D\in Z^1_\Ad(G,\cg\tens \cg)$ and extending the
latter to an invariant bivector defines a Poisson bracket. Working
backwards one has axioms of a Poisson-Lie group and its tangent
space at the identity is a Lie bialgebra. The notions are due to
Drinfeld. An introduction is in \cite{Ma:book}.

Now, for a Hopf algebra $H$ we have the following notions of
covariance of a differential structure, which we shall aim to
semiclassicalise. Note that a Hopf algebra acts and coacts on
itself from both sides via the product and coproduct respectively.
Here the coproduct is a map $\Delta:H\to H\tens H$. We will often
use the Sweedler notation $\Delta a=a\o\tens a\t$,
$(\id\tens\Delta)\Delta a=a\o\tens a\t\tens a\thr$, etc., for it.

\begin{defin} A differential calculus $(\Omega^1,\extd)$ on a Hopf
  algebra
$H$ is called {\em left covariant} if

\smallskip
1) \quad
  $\Omega^1$ is a left comodule for $H$ with
coaction $\Delta_L:\Omega^1\to H\tens\Omega^1$.

\smallskip
2)\quad $\Delta_L:\Omega^1\to H\tens \Omega^1$ is a bimodule map
(with the tensor bimodule structure on $ H\tens \Omega^1$).

\smallskip
3)\quad The map $\extd :H\to \Omega^1$ is a comodule map.

Similarly for a right-covariant calculus with structure map
$\Delta_R:\Omega^1\to\Omega^1\tens H$. A calculus is called
bicovariant if it is both left and right covariant and

\smallskip
4)\quad The left and right coactions on $\Omega^1$ commute (a
bicomodule). \end{defin}

Note that condition 3) here fully determines the left coaction
(similarly the right) if they exist, since every element of
$\Omega^1$ is a linear combination of ones of the form $a\extd b$,
so they are induced canonically from the `group translation' or
coproduct in the invariant case.

\subsection{Semiclassical left covariance condition}

The functions on a Poisson Lie group $G$ typically deform to a
noncommutative Hopf algebra $H$, with the commutator of two
functions being given by the Poisson bracket to first order in
$\hbar$. We suppose that the differential calculus on $H$ is a
deformation of the standard (commutative) differential calculus in
the sense that we have discussed earlier. We now ask what left,
right (and later bi) covariance of the differential calculus means
in terms of the connection $\gamma$. To this end we use the
following lemma for left (similarly right) covariant calculi on a
Hopf algebra $H$:

\begin{lemma}\label{ppoo}
     Write $\Delta_L(\tau)=\tau\bo\tens\tau\bt$ and
$\Delta_R(\tau)\,=\,\tau\bz\tens \tau\bo$.  Then for $a\in H$ and
$\tau\in\Omega^1$,
\begin{eqnarray*}\label{binv}
(({\rm id}\tens[\ ,\ ])\circ\Delta_L-\Delta_L\circ[\ ,\ ])(a\tens
\tau) &=&-\,[a\o , \tau\bo]\tens  \tau\bt a\t \,\in\,
H\tens\Omega^1\ ,\cr (([\ ,\ ]\tens{\rm id})\circ\Delta_R\,-\,
\Delta_R\circ[\ ,\ ] )(a\tens \tau) &=& \tau\bz a\o \tens
[\tau\bo,a\t ]\,\in\, \Omega^1\tens H\ . \end{eqnarray*}
\end{lemma}

\proof These results follow from the following equations:
\begin{eqnarray*}
\Delta_L[a,\tau] &=& \Delta_L(a\tau-\tau a)\,=\,a\o \tau\bo\tens
a\t \tau\bt\,-\,\tau\bo a\o \tens \tau\bt a\t \ , \cr ({\rm
id}\tens[\ ,\ ])\Delta_L(a\tens \tau) &=&({\rm id}\tens[\ ,\
])(a\o  \tau\bo\tens a\t \tens
\tau\bt)\\
&=&\,a\o  \tau\bo\tens (a\t  \tau\bt- \tau\bt a\t ) \cr
\Delta_R[a,\tau] &=& \Delta_R(a\tau-\tau a)\,=\,a\o \tau\bz\tens
a\t \tau\bo\,-\,\tau\bz a\o \tens \tau\bo a\t \ , \cr ([\ ,\
]\tens{\rm id})\Delta_R(a\tens \tau) &=&([\ ,\ ]\tens{\rm
id})(a\o \tens \tau\bz\tens a\t  \tau\bo)\\
&=&\, (a\o \tau\bz-\tau\bz a\o )\tens a\t  \tau\bo\ .\square
\end{eqnarray*}

For a Lie group $G$, there is a left $G$-action on functions on
$G$ given by $(g\la a)(h)=a(g^{-1}h)$. The left multiplication map
$L_g:G\to G$ has a derivative $L_{g*}:T_hG\to T_{gh}G$ at all
$h\in G$, and the dual of this map is $L_g^*:T_{gh}^*G\to T_h^*G$.
There is a left action on $\tau\in T^*G$ given by $(g\la
\tau)(h)=L_{g^{-1}}^*(\tau(g^{-1}h))$ for all $h$. In terms of the
algebra of functions on $G$, we have (in an appropriate setting) a
left coaction $(\Delta_L a)(g,h)\,=\,a(gh)=(g^{-1}\la a)(h)$, and
a left coaction on the 1-forms $(\Delta_L
\tau)(g,h)\,=L_{g}^*(\tau(gh))=(g^{-1}\la\tau)(h)$. In our case we
work with $C^\infty(G)$ rather than algebraically, hence we work
directly with actions and group multiplication rather than such
coactions and coproduct. When working abstractly on the group (or
any other) manifold we employ the notations
\[ v(a)(g)=D_{(g;v(g))}a=a'(g;v(g))={\extd \over\extd t}|_0 a(g(t))\]
for the action of a vector field $v$, where $g(t)$ is a curve with
$g(0)=g$ and tangent $v(g)$ there. Also, if $v\in\cg=T_eG$ we
write the associated left invariant vector field $L_{*v}$ with
values $L_{g*}v\in T_gG$, or simply $gv$ in a shorthand. We
similarly have right multiplication $R_g:G\to G$, etc., and right
invariant vector field $R_{*v}$ generated by $v\in \cg$.

To interpret the left invariance formula in Lemma~\ref{ppoo} for a
Poisson Lie group $G$, we remember that, up to $O(\hbar)$, we have
same coaction formulae as in the Lie group case. We also need to
trivialise the cotangent bundle by a map $T^*G\to G\times \cg^*$
defined by $(g;\xi)\mapsto (g,L^*_g\xi)$ for $\xi\in T_g^*G$.
Given a section $\tau$ of $T^*G$ we define a section $\tilde \tau$
of the trivial $\cg^*$ bundle (a function on $G$ with values in
$\cg^*$) by \eqn{conv}{ \tilde\tau(g)=L^*_g
(\tau(g))=(g^{-1}\la\tau)(e)} where $e$ is the group identity.

\begin{lemma}\label{llinv} The semiclassicalisation $\gamma$ of a
left-invariant calculus on $G$ obeys
\begin{eqnarray*}
\tilde\gamma(a,\tilde\tau)(gh)\,-\,\tilde\gamma(g^{-1}\la a,
\widetilde{g^{-1}\la\tau})(h)\,=\, a'(gh;\omega\so (g)h)\
\tilde\tau'(gh;\omega\st (g)h)\  \end{eqnarray*} as the lowest
order part of the condition in Lemma~\ref{ppoo}.
\end{lemma}
\proof To calculate $[a\o , \tau\bo]\tens \tau\bt a\t $ we use the
Poisson bracket with bivector $\omega=\omega\so \tens\omega\st $
whereby $\{a,b\}=D_{(g;\omega\so(g))}(a)D_{(g;\omega\st(g))}(b)$
for functions $a,b$. Then
\begin{eqnarray*}
([a\o , \tau\bo]\tens  \tau\bt a\t )(g,h) &=&
\hbar\,(D_{(g;\omega\so (g))} a(gh))\, (D_{(g;\omega\st (g))}
L^*_g (\tau(gh)))+ O(\hbar^2) \cr
  \Delta_L[a,\tau](g,h) &=&
\hbar\,L^*_g(\gamma(a,\tau)(gh))\,+\, O(\hbar^2)\\(\id\tens[\ ,\
])\Delta_L(a\tens\tau)(g,h)&=& \hbar\, \gamma(g^{-1}\la
a,g^{-1}\la\tau)(h)\,+\,O(\hbar^2)\ .
\end{eqnarray*} Now we can evaluate the first half of Lemma~\ref{ppoo} at
$(g,h)$ to order $\hbar$ and apply $L^*_h$ to both sides to get
\begin{eqnarray*}&&\nquad L^*_{gh}(\gamma(a,\tau)(gh))\,-\,
L^*_h(\gamma(g^{-1}\la a,g^{-1}\la\tau)(h))\\&&=\,
(D_{(g;\omega\so (g))} a(gh))\, (D_{(g;\omega\st (g))}
L^*_{gh}(\tau(gh)))\ .
\end{eqnarray*}
This is \begin{eqnarray*} &&\nquad
((gh)^{-1}\la\gamma(a,\tau))(e)-(h^{-1}\la\gamma(g^{-1}\la
a,g^{-1}\la\tau))(e)\\
&&=(D_{(g;\omega\so (g))} a(gh))\, (D_{(g;\omega\st (g))}
((gh)^{-1}\la\tau)(e)))\end{eqnarray*} which is the equation
stated. Note that by definition
$\tilde\gamma(a,\tilde\tau)=\widetilde{\gamma(a,\tau)}$ so that
everything in the $\tilde{\ }$ notation is referred to functions
on $G$. \eproof

Conversely, if we have only a preconnection
as semiclassical datum for a calculus, we say that it is 
left-invariant if the condition in the lemma holds. We now give
a different characterisation.

\begin{defin}\label{connen}
On the trivial bundle $G\times\cg^*\to G$ we define the
bilinear-valued function $\Xi:G\times\cg^*\times\cg^*\to\cg^*$ by
\begin{eqnarray*}
\tilde\gamma(a,s)(g)\,=\,a'(g;\omega\so (g))\, s'(g;\omega\st (g))
\,+\,\Xi(g;\hat L_a(g),s(g))\,=\,\{a,s\}(g)\,+\,\Xi(g,\hat
L_a(g),s(g)) \ , \end{eqnarray*} where $a\in C^\infty(G)$, $s\in
C^\infty(G,\cg^*)$ and $\hat L_a:G\to\cg^*$ is the left invariant
derivative $\hat L_a(g)(v)=a'(g;gv)$. \end{defin} In the second
expression we have extended the notation for the Poisson bracket
to include $\cg^*$ valued functions, and also note that $\hat
L_a\,=\,\widetilde{\extd a}$.

\begin{propos} A preconnection $\gamma$ is left invariant
  if and only if $\Xi:G\times \cg^*\times \cg^*\to \cg^*$
is independent of the $G$ variable. \end{propos} \proof To apply
the condition in Lemma~\ref{llinv}, we need to note that
$(\widetilde{g^{-1}\la \tau})(h)=\tilde\tau(gh)$, and then
\begin{eqnarray*}
&&a'(gh;\omega\so (gh))\, \tilde\tau'(gh;\omega\st (gh)) +\Xi(gh,
\hat L_a(gh),\tilde\tau(gh))-a'(gh;g\,\omega\so (h))\,
\tilde\tau'(gh;g\,\omega\st (h)) \cr && - \Xi(h,\hat L_{g^{-1}\la
a}(h),\tilde\tau(gh))\,=\, a'(gh;\omega\so
(g)h)\,\tilde\tau'(gh;\omega\st (g)h)\ .
\end{eqnarray*} In $TG\tens TG$ we have
$\omega(gh)-g\,\omega(h)-\omega(g)\,h=0$, so the last equation
reduces to
\begin{eqnarray*}
\Xi(gh, \hat L_a(gh),\tilde\tau(gh))\,=\, \Xi(h,\hat L_{g^{-1}\la
a}(h), \tilde\tau(gh))\ . \end{eqnarray*} A bit of calculation
shows that $\hat L_{g^{-1}\la a}(h)=\hat L_a(gh)$, so we are left
with
\begin{eqnarray*}
\Xi(gh, \phi,\psi)\,=\, \Xi(h,\phi,\psi)\quad
\forall\phi,\,\psi\in \cg^* \quad \forall g,\,h\in G\ .
\end{eqnarray*} We deduce that $\Xi(h,\phi,\psi)$ is independent of
$h\in G$. \eproof

\begin{propos}\label{uutt}
A left invariant preconnection $\tilde\gamma$ corresponding
to $\Xi:\cg^*\tens\cg^*\to\cg^*$ is compatible in the sense of
(\ref{eqq2}) if and only if
\begin{eqnarray*}
     \Xi(\phi,\psi)\,-\,\Xi(\psi,\phi)\,=\,[\phi,\psi]_{\cg^*} \ .
\end{eqnarray*} \end{propos} \proof The left invariant
trivialisation of $T^*G$ gives the following form of (\ref{eqq2}):
\[
\tilde\gamma(x,\hat L_y)\,-\,\tilde\gamma(y,\hat L_x)\,=\,\hat
L_{\{x,y\}}\ .
\]
This can be rearranged to give
\begin{eqnarray}\label{wasyy}
\Xi(\hat L_x,\hat L_y)\,-\,\Xi(\hat L_y,\hat L_x)\,=\,\hat
L_{\{x,y\}}\,-\, \{x,\hat L_{y}\}\,+\,\{y,\hat L_{x}\}\ .
\end{eqnarray}
We only have to evaluate this equation at the identity $e\in G$.
All Poisson brackets evaluated at the identity vanish as
$\omega(e)=0$, so to find the right hand side of (\ref{wasyy}) we
only need
\begin{eqnarray*}
\hat L_{\{x,y\}}(g)(v) &=& D_{(g;gv)}\,\{x,y\}\ ,\cr \hat
L_{\{x,y\}}(e)(v) &=& x'(e;{\omega\so }'(e;v))\, y'(e;\omega\st
(e))\,+\, x'(e;\omega\so (e))\,y'(e;{\omega\st }'(e;v))\ .
\end{eqnarray*} Setting $\phi(v)=x'(e;v)$ and $\psi(v)=y'(e;v)$,
we see that evaluating (\ref{wasyy}) at $e$ gives
\begin{eqnarray*}
\Xi(\phi,\psi)(v)\,-\,\Xi(\psi,\phi)(v)\,=\,(\phi\tens\psi)(\omega'(e;v))
\,=\,(\phi\tens\psi)(\delta(v))\,=\,[\phi,\psi]_{\cg^*}(v) \
.\quad\square \end{eqnarray*}

\subsection{Left covariance in the quasitriangular case}

The quasitriangular case is the most important because the
standard q-deformation quantum groups quantise Poisson-Lie groups
of this type. In this case there are some useful simplifications
which will be needed in the next section. We recall that a Lie
bialgebra $\cg$ is called quasitriangular if the Lie cobracket
$\delta$ is of the form $\delta v=\ad_v(r)$ for all $v\in \cg$,
where $r=r\so \tens r\st \in\cg\tens\cg$ obeys the Classical
Yang-Baxter equations
\[ [[r,r]]\equiv
[r_{12},r_{13}]+[r_{12},r_{23}]+[r_{13},r_{23}]=0.\] The
expression here is the Schouten bracket and the numbers refer to
the position in $\cg\tens\cg\tens\cg$. In this case
$\omega(g)\,=\,g\,r\,-\,r\,g$ for all $g\in G$. Set $r_{\pm}=(r\pm
r_{21})/2$, where $r_{21}=r\st \tens r\so $. Then $r_+$ is
necessarily $\ad$-invariant, as is the element
\[ n\equiv[[r_-,r_-]]=-[[r_+,r_+]]=[r_{+12},r_{+23}].\] A 
quasitriangular Lie bialgebra is called triangular if $r_+=0$.

\begin{cor}\label{uuttii} If $\cg$ is quasitriangular, a left
invariant preconnection $\tilde\gamma$ is compatible if and
only if the map $\hat\Xi:\cg^*\tens\cg^*\to \cg^*$ defined by
\[\hat\Xi(\phi,\psi) =\Xi(\phi,\psi)\,-\,
\phi({r_-}\so )\,\ad^{*}_{{r_-}\st }\psi
     \]
is symmetric. Here $\ad^*_v\psi(w)=\psi([w,v])$.
\end{cor} \proof  Because $r_+$ is ad-invariant, the
result of Proposition~\ref{uutt} can be written as
\begin{eqnarray*}
\Xi(\phi,\psi)(v)\,-\,\Xi(\psi,\phi)(v)\,=\,(\phi\tens\psi)([v,{r_-}\so
]\tens{r_-}\st + {r_-}\so \tens [v,{r_-}\st ]) \ .
\end{eqnarray*} Then we have the
following, which is zero by the antisymmetry of $r_-$:
\begin{eqnarray*}
\hat\Xi(\phi,\psi)(v)\,-\,\hat\Xi(\psi,\phi)(v) &=&
\phi([v,{r_-}\so ])\,\psi({r_-}\st )\,+\, \phi({r_-}\so )\,\psi(
[v,{r_-}\st ]) \cr && -\, \phi({r_-}\so )\,\psi([v,{r_-}\st
])\,+\, \psi({r_-}\so )\,\phi([v,{r_-}\st ])
  \ .\quad\square
\end{eqnarray*}

\begin{propos}\label{propj1} In the quasi-triangular case,
the first super Jacobi identity is equivalent to the following
equation for all $\phi,\psi,\zeta\in\cg^{*}$,
     \begin{eqnarray*}
     \Xi(\phi,\Xi(\psi,\zeta))\,-\,\Xi(\psi,\Xi(\phi,\zeta)) &=&
     \phi({r_-}\so )\, \Xi(\ad^{*}_{{r_-}\st }\psi,\zeta) \, +\,
     \psi({r_-}\st )\,\Xi(\ad^{*}_{{r_-}\so }\phi,\zeta)\ .
\end{eqnarray*} \end{propos} \proof The first super Jacobi
identity reduces to showing that
\begin{eqnarray}\label{uui9}
     \tilde\gamma(x,\tilde\gamma(y,s))\,-\,
     \tilde\gamma(y,\tilde\gamma(x,s)) &=& \tilde\gamma(\{x,y\},s)\ .
\end{eqnarray} If we use the Poisson bracket notation, we can
write
\begin{eqnarray}\label{hsighh}
     \tilde\gamma(x,\tilde\gamma(y,s))(g) &=& \{x,\{y,s\}\}(g)\,
     +\,\{x,\Xi(\hat L_y(g),s(g))\} \cr
&& +\, \Xi(\hat L_x(g),\{y,s\}(g))\,+\,\Xi(\hat L_x(g),\Xi(\hat
L_y(g),s(g)))\ . \end{eqnarray} Then the Jacobi identity for the
Poisson bracket shows that all the double Poisson bracket terms in
(\ref{uui9}) cancel, so we get
\begin{eqnarray*}
\Xi(\hat L_{\{x,y\}}(g),s(g)) &=& \{x,\Xi(\hat L_y(g),s(g))\} \,
+\, \Xi(\hat L_x(g),\{y,s\}(g))\\
&&+\,\Xi(\hat L_x(g),\Xi(\hat L_y(g),s(g))) - \{y,\Xi(\hat
L_x(g),s(g))\} \\&& -\, \Xi(\hat L_y(g),\{x,s\}(g))\,-\,\Xi(\hat
L_y(g),\Xi(\hat L_x(g),s(g)))\ .
\end{eqnarray*} Using the equation $\{x,\Xi(\hat
L_y(g),s(g))\}=\Xi(\{x,\hat L_y\}(g),s(g))+\Xi(\hat
L_y(g),\{x,s\}(g))$, we find
\begin{eqnarray*}
\Xi(\hat L_{\{x,y\}}(g),s(g)) &=& \Xi(\{x,\hat
L_y\}(g),s(g))\,+\,\Xi(\hat L_x(g),\Xi(\hat L_y(g),s(g))) \cr &&
-\, \Xi(\{y,\hat L_x\}(g),s(g)) \,-\,\Xi(\hat L_y(g),\Xi(\hat
L_x(g),s(g)))\ . \end{eqnarray*} Now calculate
\begin{eqnarray*}
\hat L_{\{x,y\}}(g)(v) &=& D_{(g;gv)}(\, x'(g;\omega\so (g))\,
y'(g;\omega\st (g))\,)\ , \cr \{x,\hat L_y\}(g)(v) &=&
x'(g;\omega\so (g))\, D_{(g;\omega\st (g))} \, y'(g;gv)\ ,
\end{eqnarray*} so using the antisymmetry of $\omega$ and the Lie
bracket of vector fields,
\begin{eqnarray*}
\hat L_{\{x,y\}}(g)(v) -\{x,\hat L_y\}(g)(v)+\{y,\hat L_x\}(g)(v)
&=& x'(g;\omega\so (g))\, y'(g;[gv,\omega\st (g)]) \cr && +
x'(g;[gv,\omega\so (g)]) \,y'(g;\omega\st (g)) .
\end{eqnarray*} Using Ad-invariance of $r_+$ we write
$\omega=g\,r_-\,-\,r_-\,g$, which gives the answer. \eproof

We end with some lemmas to be used later.

\begin{lemma}\label{thhh}
     In terms of $\hat\Xi$, the first super Jacobi identity is
  equivalent to
     the vanishing of the following expression
for all $\phi,\psi,\zeta\in\cg^*$, where $n=[\![r_-,r_-]\!]$:
\begin{eqnarray*}
     &&\hat \Xi(\phi,\hat\Xi(\psi,\zeta))\,-\, \hat
     \Xi(\psi,\hat\Xi(\phi,\zeta))
     \,+\,\phi(n\so )\,\psi(n\st )\,\ad^{*}_{n\sth }\,\zeta\cr &&
+\,\phi({r_-}\so )\,(\,\ad^{*}_{{r_-}\st }\,\hat\Xi(\psi,\zeta)
\,-\, \hat
     \Xi(\psi,\ad^{*}_{{r_-}\st }\,\zeta)\,-\,
     \hat\Xi(\ad^{*}_{{r_-}\st }\psi,\zeta)\,)\cr
&&-\,\psi({r_-}\so )\,(\,\ad^{*}_{{r_-}\st
}\,\hat\Xi(\phi,\zeta)\,-\, \hat
     \Xi(\phi,\ad^{*}_{{r_-}\st }\,\zeta)\,-\,
     \hat\Xi(\ad^{*}_{{r_-}\st }\phi,\zeta)\,)\ .
\end{eqnarray*} \end{lemma} \proof Let $t$ be an independent copy
of $r$. By definition of $\hat\Xi$,
\begin{eqnarray*}
     \Xi(\phi,\Xi(\psi,\zeta)) &=& \hat \Xi(\phi,\hat\Xi(\psi,\zeta))\,+\,
     \psi({r_-}\so )\,\hat \Xi(\phi,\ad^{*}_{{r_-}\st }\,\zeta) \\
     &&+\phi({r_-}\so )\,\ad^{*}_{{r_-}\st }\,\hat\Xi(\psi,\zeta)
     +\phi({r_-}\so )\,\psi({t_-}\so )\,\ad^{*}_{{r_-}\st }\,
     \ad^{*}_{{t_-}\st }\,\zeta,\cr \phi({r_-}\so )\,\Xi(
     \ad^{*}_{{r_-}\st }\psi,\zeta) &=&
     \phi({r_-}\so )\,\hat\Xi(\ad^{*}_{{r_-}\st }\psi,\zeta)\,+\,
     \phi({r_-}\so )\,(\ad^{*}_{{r_-}\st }\psi)({t_-}\so )\,
     \ad^{*}_{{t_-}\st }\zeta\ .
\end{eqnarray*} By Proposition~\ref{propj1}, we get the vanishing of
\begin{eqnarray*}
     &&\hat \Xi(\phi,\hat\Xi(\psi,\zeta))\,-\, \hat
     \Xi(\psi,\hat\Xi(\phi,\zeta)) \,+\, \psi({r_-}\so )\,\hat
     \Xi(\phi,\ad^{*}_{{r_-}\st }\,\zeta) \,-\, \phi({r_-}\so )\,\hat
     \Xi(\psi,\ad^{*}_{{r_-}\st }\,\zeta) \cr &&+\,
     \phi({r_-}\so )\,\ad^{*}_{{r_-}\st }\,\hat\Xi(\psi,\zeta)\,-\,
     \psi({r_-}\so )\,\ad^{*}_{{r_-}\st }\,\hat\Xi(\phi,\zeta)\,-\,
     \phi({r_-}\so )\,\hat\Xi(\ad^{*}_{{r_-}\st }\psi,\zeta) \cr && +\,
     \psi({r_-}\so )\,\hat\Xi(\ad^{*}_{{r_-}\st }\phi,\zeta)\,+\,
     \phi(n\so )\,\psi(n\st )\,\ad^{*}_{n\sth }\,\zeta\ ,
\end{eqnarray*} and this can be rearranged to give the
answer. \eproof

\begin{lemma}\label{upup} Given $\hat\Xi$, we recover $\gamma$ by
\[ \gamma(x,\extd y)(g)=  x'(g;\omega\so (g)) \extd
y'(g;\omega\st (g)) +L^*_{g^{-1}}\hat\Xi(\hat L_{x}(g),\hat
L_{y}(g)).\]
\end{lemma} \proof Evaluating Definition~\ref{connen} against
$v\in\cg$ we have
\begin{eqnarray*}
     \tilde\gamma(x,\hat L_{y})(g)(v) &=&\{x,\hat L_y\}(g)(v)+\Xi(\hat
L_{x}(g),\hat
     L_{y}(g))(v)\\
     &=&x'(g;\omega\so (g))\, D_{(g;gv)}\,y'(g;\omega\st (g))\,+\,
     x'(g;g\,r_-\so )\,y'(g;g[r_-\st ,v])\cr &&+\,
     \Xi(\hat L_{x}(g),\hat L_{y}(g))(v)  \cr
&=&
     x'(g;\omega\so (g))\, D_{(g;gv)}\,y'(g;\omega\st (g))\,+\,
     \hat\Xi(L_{x}(g),L_{y}(g))(v)\ .
\end{eqnarray*}  We then convert back from $\tilde\gamma$ to
$\gamma$ via (\ref{conv}). \eproof

For the remaining lemmas, as well as elsewhere, it is convenient
to switch to an alternative notation for the bilinear functions
$\Xi$ and $\hat\Xi$ from $\cg^*\tens\cg^*$ to $\cg^*$. We can
consider them as linear functions from $\cg$ to $\cg\tens \cg$
with notation
\[
\Xi(\phi,\psi)(v)\,=\,(\phi\tens\psi)(\Xi(v))\ ,\quad
\hat\Xi(\phi,\psi)(v)\,=\,(\phi\tens\psi)(\hat\Xi(v))
\]
where $\Xi(v)\,=\,\Xi\so(v)\tens \Xi\st(v)$ and
$\hat\Xi(v)\,=\,\hat\Xi\so(v)\tens \hat\Xi\st(v)$. Then we see
that $\hat\Xi(v)=\Xi(v)+{r_-}\so \tens [{r_-}\st ,v]$.

\begin{lemma}\label{nameless} The commutator between a 2-form
$\tilde\tau$ (expressed in the left trivialisation) and a function
$y$ is
     given by
\[
\tilde\gamma(y,\tilde\tau)(g)(v,w)\,=\,y'(g;\omega\so
(g))\,\tilde\tau'(g;\omega\st (g))\,+\, (\hat
L_{y}(g)\tens\tilde\tau(g))(\Xi(v)\tens w-\Xi(w)\tens v)\ .
\]
\end{lemma} \proof Setting
$\tau=\xi\wedge\eta$ for two 1-forms $\xi$ and $\eta$,
\begin{eqnarray*}
\tilde\gamma(y,\tau) &=& \tilde\gamma(y,\xi)\wedge\eta\,+\,
\xi\wedge\tilde\gamma(y,\eta)\cr \tilde\gamma(y,\tau)(g)&=&
y'(g;\omega\so (g))\,\Big(\xi'(g;\omega\st (g)) \wedge\eta(g)+
\xi(g)\wedge\eta'(g;\omega\st (g))\Big)\cr &&+\, \Xi(\hat
L_{y}(g),\xi(g))\wedge\eta(g)\,+\, \xi(g)\wedge \Xi(\hat
L_{y}(g),\eta(g))\ .
\end{eqnarray*} Applying this to $(v,w)$ gives
\begin{eqnarray*}
\tilde\gamma(y,\tau)(g)(v,w) &=& y'(g;\omega\so
(g))\,\tau'(g;\omega\st (g))(v,w)
+ \hat L_{y}(g)(\Xi\so (v))\,\xi(g)(\Xi\st (v))\,\eta(g)(w)\\
&&\kern -20pt -\, \hat L_{y}(g)(\Xi\so (w))\,\xi(g)(\Xi\st
(w))\,\eta(g)(v)
+ \xi(g)(v)\,\hat L_{y}(g)(\Xi\so (w))\, \eta(g)(\Xi\st (w))\\
&&\kern -20pt -\, \xi(g)(w)\,\hat L_{y}(g)(\Xi\so (v))\,
\eta(g)(\Xi\st (v))\ ,
\end{eqnarray*} which can be rearranged to give the
answer. \eproof

\begin{lemma}\label{string}
For a 1-form $\xi$ expressed in the left trivialisation, the
$\extd $ operation is given by
\[
\widetilde{\extd \xi}(g)(v,w) \,=\,
D_{(g;gv)}\tilde\xi(g)(w)\,-\,D_{(g;gw)}\tilde\xi(g)(v)\,+\,
\tilde\xi(g)([w,v])\ .
\]
Also we have
\begin{eqnarray*}
&&\quad\nquad \extd
\,D_{(g;X)}\,\tilde\xi(v,w)\,-\,D_{(g;X)}\,\extd \tilde\xi(v,w)\\
&&=\, \tilde\xi'(g;D_{(g;gv)}\,X(g)-Xv)(w)\,-\,
\tilde\xi'(g;D_{(g;gw)}\,X(g)-Xw)(v)\ .
\end{eqnarray*}

\end{lemma} \proof By definition $\extd
\xi(g)(gv,gw)=\xi'(g;gv)(gw)-\xi'(g;gw)(gv)$. However
\[
D_{(g;gv)}\tilde\xi(g)(w)\,=\,D_{(g;gv)}(\xi(g)(gw))\,
=\,\xi'(g;gv)(gw)\,+\,\xi(g)(gvw)\ .\quad\square
\]

\begin{lemma}\label{rope} The commutator between a 2-form
$\widetilde{\extd \xi}$ (expressed in the left trivialisation) and
a function $y$ is
     given by
\begin{eqnarray*}
\tilde\gamma(x,\extd {\tilde\xi})(v,w) &=& x'(g;\omega\so
(g))\,\extd \,D_{(g;\omega\st (g))}\,{\tilde\xi}(v,w)
  \,+\,\hat L_{x}({r_-}\so )\,\extd
(\ad^{*}_{{r_-}\st }({\tilde\xi}))(g)(v,w) \cr && +\,
  (\hat L_{x}\tens \extd \,{\tilde\xi})(\hat\Xi(v)\tens w-
\hat\Xi(w)\tens v)\ .
\end{eqnarray*} \end{lemma} \proof We begin by
\begin{eqnarray*}
\tilde\gamma(x,\extd {\tilde\xi})(v,w) &=& x'(g;\omega\so
(g))\,D_{(g;\omega\st (g))}\,\extd \,{\tilde\xi}(v,w)\,+\, (\hat
L_{x}\tens \extd \,{\tilde\xi})(\Xi(v)\tens w-\Xi(w)\tens v) \cr
  &=& x'(g;\omega\so (g))\,\extd \,D_{(g;\omega\st (g))}\,
{\tilde\xi}(v,w)
  \,-\,x'(g;\omega\so (g))\,{\tilde\xi}'(g;D_{(g;gv)}\,
\omega\st (g)\\
  &&-\omega\st (g)v)(w)
+\,x'(g;\omega\so (g))\,{\tilde\xi}'(g;D_{(g;gw)}\,\omega\st
(g)-\omega\st (g)w)(v)
  \\ &&+\,
(\hat L_{x}\tens \extd \,{\tilde\xi})(\Xi(v)\tens w-\Xi(w)\tens v)
\cr
  &=& x'(g;\omega\so (g))\,\extd \,D_{(g;\omega\st (g))}\,
{\tilde\xi}(v,w)
  \,-\,\hat L_{x}({r_-}\so )\,{\tilde\xi}'(g;g[v,{r_-}\st ])(w) \cr &&
  +\,\hat L_{x}({r_-}\so )\,{\tilde\xi}'(g;g[w,{r_-}\st ])(v) \,+\,
(\hat L_{x}\tens
  \extd \,{\tilde\xi})(\Xi(v)\tens w-\Xi(w)\tens v) \cr
   &=& x'(g;\omega\so (g))\,\extd \,D_{(g;\omega\st (g))}\,
{\tilde\xi}(v,w)
  \,-\,\hat L_{x}({r_-}\so )\,{\tilde\xi}'(g;g[v,{r_-}\st ])(w) \cr
  &&
  +\,\hat L_{x}({r_-}\so )\,{\tilde\xi}'(g;g[w,{r_-}\st ])(v) \,+\,
(\hat L_{x}\tens
  \extd \,{\tilde\xi})(\hat\Xi(v)\tens w-\hat\Xi(w)\tens v) \cr &&
  +\, \hat L_{x}({r_-}\so )\,
  \extd {\tilde\xi}([v,{r_-}\st ],w)\,-\, \hat L_{x}({r_-}\so )\,
  \extd {\tilde\xi}([w,{r_-}\st ],v) \cr
  &=& x'(g;\omega\so (g))\,\extd \,D_{(g;\omega\st (g))}\,{\tilde\xi}(v,w)
  \,-\,\hat L_{x}({r_-}\so )\,{\tilde\xi}'(g;gw)([v,{r_-}\st ]) \cr &&
  +\,\hat L_{x}({r_-}\so )\,{\tilde\xi}'(g;gv)([w,{r_-}\st ]) \,+\, (\hat
L_{x}\tens
  \extd \,{\tilde\xi})(\hat\Xi(v)\tens w-\hat\Xi(w)\tens v)\cr
  &&-\, \hat
L_{x}({r_-}\so )\,{\tilde\xi}(g)([[v,{r_-}\st ],w]-[[w,{r_-}\st
],v])\cr
  &=&
  x'(g;\omega\so (g))\,\extd \,D_{(g;\omega\st (g))}\,{\tilde\xi}(v,w)
  \,+\,\hat L_{x}({r_-}\so )\,\extd
(\ad^{*}_{{r_-}\st }({\tilde\xi}))(g)(v,w) \cr && +\,
  (\hat L_{x}\tens \extd \,{\tilde\xi})(\hat\Xi(v)\tens w-\hat\Xi(w)
\tens v)\ .\quad\square
\end{eqnarray*}

\subsection{Example: quantum $SU_2$ with its left-covariant calculus}

Here we pause to compute a quantum group example, namely the
3-dimensional calculus on $C_\hbar(SU_2)$ found in \cite{Wor:dif}.
As with our previous example, we do not need the full $C^*$ theory
but just a subalgebra generated by a matrix of generators
$\begin{pmatrix}a&b\\
c & d\end{pmatrix}$ which in the classical limit become the matrix
element coordinate functions on $SU_2$. We use the standard
relations
\[ba=qab,\quad ca=qac,\quad db=qbd,\quad dc=qcd,\quad cb=bc,\quad
ad=1+q^{-1}bc,\quad da=1+qbc.\] For our purposes $q=e^{\hbar\over
2 }$ rather than $q$ real as needed for completion to a $C^*$
algebra. For the rest of the Hopf algebra structure we refer to
any textbook on quantum groups.  The 3-d calculus has a basis of
left-invariant 1-forms
\[ \tau^+=d\extd b-q b\extd d,\quad \tau^-=q a\extd
c-c\extd a,\quad \tau^3=d\extd a-q b\extd c\] with commutation
relations
\[ \tau^\pm \begin{pmatrix}a&b\\ c &
d\end{pmatrix}=\begin{pmatrix}qa&q^{-1}b\\ qc &
q^{-1}d\end{pmatrix}\tau^\pm,\quad \tau^3
\begin{pmatrix}a&b\\ c &
d\end{pmatrix}=\begin{pmatrix}q^2 a&q^{-2}b\\ q^2c &
q^{-2}d\end{pmatrix}\tau^3\ .\] One has
\[ \extd a=a\tau^3+q^{-1}b\tau^-,\quad \extd
b=a\tau^+-q^{-2}b\tau^3,\quad \extd c=c \tau^3+q^{-1}d\tau^-,\quad
\extd d=c\tau^+-q^{-2}d\tau^3.\]

The Poisson bracket can be read off from the algebra relations as
generated by the nonzero values
\[ \{a,b\}=-{1\over 2}ab,\quad \{a,c\}=-{1\over 2}ac,\quad
\{a,d\}=-bc,\quad \{b,d\}=-{1\over 2}bd,\quad \{c,d\}=-{1\over
2}cd.\] Classically the $\tau^\pm,\tau^3$ are dual to the left
invariant vector fields $\del_\pm,\del_3$ generated by the Lie
algebra $su_2$ and
\begin{eqnarray}\label{deriv91}
\extd x= (\del_ix)\tau^i\ ,
\end{eqnarray}
  where the sum is
over $i=\pm,3$. Explicitly here, \[ \del_+\begin{pmatrix}a&b\\
c & d\end{pmatrix}=  \begin{pmatrix}0&a\\
0 & c\end{pmatrix},\quad \del_-\begin{pmatrix}a&b\\
c & d\end{pmatrix}=\begin{pmatrix}b&0\\
d & 0\end{pmatrix},\quad \del_3\begin{pmatrix}a&b\\
c & d\end{pmatrix}= \begin{pmatrix}a&-b\\
c & -d\end{pmatrix}\] computed by right multiplication in the
usual representation of the Chevalley basis of $su_2$. Meanwhile,
from the bimodule commutation relations we read off
\[ \gamma(\begin{pmatrix}a&b\\ c &
d\end{pmatrix},\tau^i)={1\over 2}\lambda_i\begin{pmatrix}a&-b\\ c
& -d\end{pmatrix}\tau^i,\quad\lambda_\pm=-1,\quad \lambda_3=-2\]
and since $\gamma(\ ,\tau^i)$ acts as derivation, comparing values
on these generators, one finds \eqn{gamma3d}{
\gamma(x,\tau^i)={1\over 2}\lambda_i(\del_3 x)\tau^i} for all
functions $x$. Hence for all $\eta=\eta_i\tau^i$,
\begin{eqnarray}\label{eq222}
  \hat\nabla_{x}(\eta)=\gamma(x,\eta)=\{x,\eta_i\}\tau^i+{1\over
2}\lambda_i(\del_3x) \eta_i\tau^i\ .
\end{eqnarray}
  From this and the compatibility
condition (\ref{eqq2}) one may verify that the curvature
$R(\hat x,\hat y)(\tau^i)=0$ and hence $R(\hat x,\hat y)=0$ as it
must since $\Omega^1$ for the quantum calculus is a bimodule. We
are not in the symplectic case so the torsion need not vanish.

We can understand this example using our above results as follows.
The standard quasitriangular structure is \[ r=e_+\tens e_-+
{1\over 4} e_3\tens e_3\] in the Chevalley basis (so that
$[e_3,e_\pm]=\pm 2 e_\pm$ and $[e_+,e_-]=e_3$) and one may check
that this reproduces the above Poisson-bracket with
$\omega=gr-rg$. This is the usual picture for $SU_2$ with its
standard Drinfeld-Sklyanin Poisson bracket.

Next, the Lie cobracket is
\[ \delta e_i={1\over 2}(e_i\tens e_3-e_3\tens e_i)\ ,\]
so
\[ \Xi(e_i)={1\over 2}\lambda_i e_3\tens e_i,\quad \lambda_\pm=-1\ .\]
We will solve the condition in Proposition~\ref{uutt} for any
$\lambda_3$, giving a compatible preconnection. In the left
trivialisation, suppose (as above) that $\tilde\tau^i\in \cg^*$
are constant and dual to the basis $e_i$. Then
\[ \tilde\gamma(a,\tilde\tau^i)=\Xi(\hat L_a(g),\tilde\tau^i)=\tilde\tau^j
a'(g;g\Xi\so(e_j))\<\Xi\st(e_j),\tilde\tau^i\>={1\over 2}
\lambda_i(\del_3a)\tilde\tau^i.\] Since $\tau^i$ are the left
invariant extensions of $\tilde\tau^i$ (the Maurer-Cartan form) we
have exactly (\ref{gamma3d}) except that for us $\lambda_3$ is
arbitrary. For this form of connection the compatibility condition
reduces to $\del_i\{x,y\}=\{\del_i x,y\}+\{x,\del_i y\}$ which
holds as it must. The curvature $R$ then vanishes by this
equation, for all $\lambda_3$. There is a similar natural form for
$\Xi$ for all simple $\cg$ with their standard Poisson-Lie group
structures according to the form of $\delta$.

The Poisson tensor is
\[
c\,b(\del_-\tens\del_+-\del_+\tens\del_-)
+\frac{b\,d}{2}(\del_+\tens\del_3-\del_3\tens\del_+)
+\frac{c\,a}{2}(\del_-\tens\del_3-\del_3\tens\del_-)\
.
\]
Now (summing over repeated indices), (\ref{deriv91}) and
(\ref{eq222}) give
\[
\hat\nabla_{ x}(\extd z)\,=\, \{x,\del_i z\}\,\tau^i + \frac12
\,\lambda_i (\del_3 x)(\del_i z)\,\tau^i\ ,
\]
so we get
\begin{eqnarray*}
\<\hat y,\hat\nabla_{ x}(\extd z)\> &=& \omega^{ji}(\del_j y)
\big(\{x,\del_i z\} + \frac12 \,\lambda_i (\del_3
x)(\del_i z)\big) \cr &=&\omega^{ji}\,\omega^{nm}(\del_j
y)(\del_n x) (\del_m\del_i z) + \frac12
\,\lambda_i\,\omega^{ji}(\del_j y) (\del_3 x) (\del_i
z)\ .
\end{eqnarray*}
Then from (\ref{tor1}) we can write the torsion as
\[
\<dz,T(\hat x,\hat y)\>\,=\, \frac12
\,\lambda_i\,\omega^{ji}(\del_j y) (\del_3 x) (\del_i
z) - \frac12 \,\lambda_i\,\omega^{ji}(\del_j x) (\del_3 y)
(\del_i z) \ ,
\]
which is non-zero.

\subsection{Bicovariant calculi on Poisson-Lie groups} We now
assume that the differential calculus is bicovariant, resulting in
important simplifications and our main results.

\begin{theorem}\label{yyuu} A left invariant preconnection
$\tilde\gamma$ given by $\Xi$ is bicovariant if, for all $g\in G$
and $\phi,\psi\in\cg^{*}$,
\[
\Xi(\phi,\psi)\,-\,\Ad^*_{g^{-1}}\, \Xi(\Ad^*_{g}\phi,\Ad
^*_{g}\psi) \,=\, \phi(g^{-1}\omega\so (g))\,
\ad^*_{g^{-1}\omega\st (g)}\,\psi\ .
\]
\end{theorem} \proof The section $t$ of the trivial $\cg^*$ bundle
given by trivialising $T^*G$ by the right action is related to the
section $s$ of the left trivialisation by
  $t(g)=(\Ad_{g^{-1}})^*s(g)$,
and the corresponding right invariant connection would be, for
some $\Psi:\cg^{*}\tens\cg^{*}\to\cg^{*}$,
\begin{eqnarray}\label{corr}
(\Ad_{g^{-1}})^*\tilde\gamma(y,s)(g)\,=\,\{y,t\} \,+\,\Psi(\hat
R_y(g),t(g))\ . \end{eqnarray} Here we have $\hat R_y(g)\,=\,(\Ad
_{g^{-1}})^*(\hat L_y(g))$, as can be seen from
\[
\hat R_y(g)(v)\,=\,y'(g;vg)\,=\,y'(g;g\,\Ad_{g^{-1}}v)\ .
\]
From now we will use the coadjoint action, $\Ad^*_g=(\Ad
_{g^{-1}})^*$. Then
\[
t'(g;v)\,=\, \Ad^*_g(\ad^*_{g^{-1}v}\,s(g)\,+\,s'(g;v))\ ,
\]
so (\ref{corr}) becomes
\[
  \tilde\gamma(y,s)(g)\nonumber=\,y'(g;\omega\so (g))\,
(\ad^*_{g^{-1}\omega\st (g)}\,s(g)\,+\,s'(g;\omega\st (g)))
\,+\,\Ad ^*_{g^{-1}}\,\Psi(\Ad^*_{g}\hat L_y(g),\Ad^*_{g}s(g))\
.\] Substituting Definition~\ref{connen} into this gives
\[
\Xi(\hat L_y(g),s(g))\,-\,\Ad^*_{g^{-1}}\, \Psi(\Ad^*_{g}\hat
L_y(g),\Ad^*_{g}s(g)) \,=\, y'(g;\omega\so (g))\,
\ad^*_{g^{-1}\omega\st (g)}\,s(g)\ .
\]
We can use the definition of $\hat L_y(g)$ to rewrite this as
\[
\Xi(\hat L_y(g),s(g))\,-\,\Ad^*_{g^{-1}}\, \Psi(\Ad^*_{g}\hat
L_y(g),\Ad^*_{g}s(g)) \,=\, \hat L_y(g)(g^{-1}\omega\so (g))\,
\ad^*_{g^{-1}\omega\st (g)}\,s(g)\ .
\]
Now we set $\phi=\hat L_y(g)\in\cg^*$ and
  $\psi=s(g)\in\cg^*$,
\[
\Xi(\phi,\psi)\,-\,\Ad^*_{g^{-1}}\,
\Psi(\Ad^*_{g}\phi,\Ad^*_{g}\psi) \,=\, \phi(g^{-1}\omega\so
(g))\, \ad^*_{g^{-1}\omega\st (g)}\,\psi\ .
\]
Setting $g=e$ now tells us that $\Xi=\Psi$. \eproof

\begin{propos}\label{hatXiad} In the case of $\cg$ quasitriangular,
a left-invariant preconnection $\tilde\gamma$ given by
$\hat\Xi$ is bicovariant if and only if
  $\hat\Xi$ is
Ad-invariant. The infinitesimal form of this is
     \[
     \hat\Xi([v,w])\,=\,[\hat\Xi\so(v),w]\tens \hat\Xi\st(v)
     \,+\,\hat\Xi\so(v)\tens [\hat\Xi\st(v),w]\ .
     \]
\end{propos} \proof Use Theorem~\ref{yyuu} and
$\omega(g)=g\,r_-\,-\,r_-g$. Ad-invariance means of course
$\hat\Xi(\phi,\psi)\,=\,\Ad^*_{g^{-1}}\,
\hat\Xi(\Ad^*_{g}\phi,\Ad^*_{g}\psi)$ for all $g\in G$, which we
then write as stated. \eproof

\begin{propos}\label{qbicj1} In the quasitriangular bicovariant case,
the first super Jacobi identity corresponds to the vanishing of
\[(\id\tens\hat\Xi)\hat\Xi(v)-(\tau\tens\id)(\id\tens\hat\Xi)\hat\Xi(v)+\,n\so
\tens n\st \tens [v,n\sth ]\] where $n=[[r_-,r_-]]$ and $\tau$ is
the usual flip map.
\end{propos} \proof We use the result of
Proposition~\ref{thhh} and $\hat\Xi$ now $\ad$-invariant. \eproof

In other words, this expression corresponds to the curvature of
the preconnection in the geometric picture of Section~3. The
second Jacobi identity is much more painful but turn out to also
be controlled by this element $n$:

\begin{propos} In the quasitriangular bicovariant case if the first
    super Jacobi identity holds, the second
     super Jacobi identity corresponds to the vanishing of
\begin{eqnarray*}&&\kern-20pt 2\,D_{(g;gw)}\big({\hat L}_{x}(e\so
(v))\big)\, \hat L_x(e\st (v))\, \hat L_x(e\sth(v)) -\,
2\,D_{(g;gv)}\big({\hat L}_{x}(e\so (w))\big)\, \hat L_x(e\st
(w))\, \hat L_x(e\sth (w))
  \\
     && +\,\hat L_x([e\so (v),w])\,\hat L_x(e\st (v))\,\hat L_x(e\sth (v))
     -\, \hat L_x([e\so (w),v])\,\hat L_x(e\st (w))\,\hat L_x(e\sth (w))\ .
\end{eqnarray*}
for all $x$ and $v,w\in\cg$, where $e(v)  \,=\, -\,n\so \tens n\st
\tens [v,n\sth ]$.
\end{propos} \proof For brevity we write Lemma~\ref{upup} as
\begin{eqnarray*}
[y,\extd x]_{\bullet} &=& \hbar\,(p_{y}^{1}\, \extd
\,p_{x}^{2}\,+\, h_{yx})\,+\,O(\hbar^{2})\ ,\cr [\extd y,\extd
x]_{\bullet} &=& \hbar\,(
  \extd \,p_{y}^{1}\wedge \extd \,p_{x}^{2}\,+\, \extd
\,h_{yx})\,+\,O(\hbar^{2})\ ,
\end{eqnarray*} where
\begin{eqnarray*}
         p_{y}^{1}(g) \,=\, y'(g;\omega\so (g))\ ,\quad
         p_{x}^{2}(g) \,=\, x'(g;\omega\st (g))\ ,\quad
         h_{yx}(g;gv) \,=\, \hat\Xi(\hat L_{y}(g),\hat L_{x}(g))(v)\ ,
     \end{eqnarray*}
and the second commutator result is just the differential of the
first. Then

\begin{eqnarray*}
     [z,[\extd y,\extd x]_{\bullet}]_{\bullet} &
      = & \hbar([z,\extd p^{1}_{y}]_{\bullet}\wedge \extd p^{2}_{x}\,+\,
     \extd p^{1}_{y}\wedge [z, \extd p^{2}_{x}]_{\bullet}\,
     +\,[z,\extd \,h_{yx}]_{\bullet})\,+\,O(\hbar^3) \cr
     &=& \hbar([z,\extd p^{1}_{y}]_{\bullet}\wedge \extd p^{2}_{x}\,+\,
     [z,\extd p^{1}_{x}]_{\bullet}\wedge \extd p^{2}_{y}\,
     +\,[z,\extd \,h_{yx}]_{\bullet})\,+\,O(\hbar^3)\ , \cr
     [\extd x,[z,\extd y]_{\bullet}]_{\bullet} &
     =& \hbar([\extd x,p^{1}_{z}\,\extd
p^{2}_{y}+h_{zy}]_{\bullet} )\,+\,O(\hbar^3)\cr
     &=& \hbar(-\,[p^{1}_{z},\extd x]_{\bullet}\wedge \extd p^{2}_{y}\,+\,
p^{1}_{z}\,
     [\extd x,\extd p^{2}_{y}]_{\bullet}\,+\,[\extd
x,h_{zy}]_{\bullet})\,+\,O(\hbar^3)\ .
\end{eqnarray*} Now we can calculate
\begin{eqnarray*}
  \hbar J_{2}(x,x,x)+O(\hbar^2) & = & \hbar^{-1}([x,[\extd x,\extd
x]_{\bullet}]_{\bullet} \,-\, [\extd x,[x,\extd
x]_{\bullet}]_{\bullet}\,
-\,[\extd x,[x,\extd x]_{\bullet}]_{\bullet}) \\
      & = & 2\,([x,\extd p^{1}_{x}]_{\bullet}
      +[p^{1}_{x},\extd x]_{\bullet})\wedge \extd p^{2}_{x}
      -\, 2\,p^{1}_{x}\,
[\extd x,\extd p^{2}_{x}]_{\bullet}\\
&&-\,2\,[\extd x,h_{xx}]_{\bullet}\,+\,[x,\extd
\,h_{xx}]_{\bullet}\ .
\end{eqnarray*} Now we can use an independent copy $q$ of $p$ to
give
\begin{eqnarray*}
     [p^{1}_{x},\extd x]_{\bullet} &=& \hbar(q_{p_{x}^{1}}^{1}\,\extd
\,q^{2}_{x}\,
     +\, h_{p_{x}^{1}x})\,+\,O(\hbar^2)\ ,\cr
     [x,\extd p^{1}_{x}]_{\bullet} &=& \hbar(q_{x}^{1}\,\extd
\,q^{2}_{p^{1}_{x}}\,
     +\, h_{p_{x}^{1}x})\,+\,O(\hbar^2)\ ,\cr
     [\extd q^{2}_{x},\extd x]_{\bullet} &=& \hbar(\extd \,p^{1}_{q^{2}_{x}}
     \wedge \extd \,p^{2}_{x}\,+\,
     \extd \,h_{q_{x}^{2}x})\,+\,O(\hbar^2)\ .
\end{eqnarray*} From this we get
\begin{eqnarray*}
   J_{2}(x,x,x)
      & = & 2\,(q_{p_{x}^{1}}^{1}\,\extd \,q^{2}_{x}\,+\,2\,
h_{p_{x}^{1}x}\,+\,
      q_{x}^{1}\,\extd \,q^{2}_{p^{1}_{x}}
      \,-\,q^{1}_{x}\,\extd \,p^{1}_{q^{2}_{x}})\wedge \extd p^{2}_{x}\cr &&
      -\,
      2\,p^{1}_{x}\, \extd \,h_{p_{x}^{2}x}\,-\,2\,[\extd
x,h_{xx}]_\bullet/\hbar\,
      +\,[x,\extd \,h_{xx}]_\bullet/\hbar\,+\,
O(\hbar)\ . \end{eqnarray*} Suppose that $p$ corresponds to
$\omega$ and $q$ to $\pi$, another copy of $\omega$.  Then
\begin{eqnarray*}q^{2}_{p^{1}_{x}}\,-\, p^{1}_{q^{2}_{x}} &=&
x'(g;D_{(g;\pi^{2}(g))}\omega\so (g)\,-\, D_{(g;\omega\so
(g))}\pi^{2}(g))\ ,\cr q_{p_{x}^{1}}^{1}\,\extd \,q^{2}_{x}\wedge
\extd p^{2}_{x} &=& p_{q_{x}^{1}}^{1}\,\extd \,p^{2}_{x}\wedge
\extd q^{2}_{x}\,=\, -\,p_{q_{x}^{1}}^{1}\,\extd \,q^{2}_{x}\wedge
\extd p^{2}_{x}\ ,\cr 2\,q_{p_{x}^{1}}^{1}\,\extd
\,q^{2}_{x}\wedge \extd p^{2}_{x} &=&
(q_{p_{x}^{1}}^{1}-p_{q_{x}^{1}}^{1})\,\extd \,q^{2}_{x}\wedge
\extd p^{2}_{x} \\
&&=\,x'(g;D_{(g;\pi^{1}(g))}\omega\so (g)\,-\, D_{(g;\omega\so
(g))}\pi^{1}(g))\,\extd \,q^{2}_{x}\wedge \extd \,p^{2}_{x}\ .
\end{eqnarray*} For $t$ and independent copy of $r$,
\begin{eqnarray*}
2\,q_{p_{x}^{1}}^{1}\,\extd \,q^{2}_{x}\wedge \extd p^{2}_{x}&=&
x'(g;g[{t_-}\so ,{r_-}\so ])\,\extd \,x'(g;g{t_-}\st )\wedge \extd
\,x'(g;g{r_-}\st )\cr &&-\, x'(g;[{t_-}\so ,{r_-}\so ]g)\,\extd
\,x'(g;{t_-}\st g)\wedge \extd \,x'(g;{r_-}\st g)\ , \cr
q_{x}^{1}\,\extd (\,q^{2}_{p^{1}_{x}} -\,p^{1}_{q^{2}_{x}})\wedge
\extd p^{2}_{x} &=& x'(g;g{t_-}\so )\,\extd \,x'(g;g[{t_-}\st
,{r_-}\so ])\wedge \extd \,x'(g;g{r_-}\st ) \cr &&-\,
x'(g;{t_-}\so g)\,\extd \,x'(g;[{t_-}\st ,{r_-}\so ]g)\wedge \extd
\,x'(g;{r_-}\st g) \cr &=&-\,x'(g;g{t_-}\so )\,\extd
\,x'(g;g{r_-}\st )\wedge \extd \,x'(g;g[{t_-}\st ,{r_-}\so ]) \cr
&&+\, x'(g;{t_-}\so g)\,\extd \,x'(g;{r_-}\st g)\wedge \extd
\,x'(g;[{t_-}\st ,{r_-}\so ]g)\cr &=&+\,x'(g;g{t_-}\so )\,\extd
\,x'(g;g{r_-}\so )\wedge \extd \,x'(g;g[{t_-}\st ,{r_-}\st ]) \cr
&&-\, x'(g;{t_-}\so g)\,\extd \,x'(g;{r_-}\so g)\wedge \extd
\,x'(g;[{t_-}\st ,{r_-}\st ]g)\ . \end{eqnarray*} Using this,
\begin{eqnarray*}
      2\,(q_{p_{x}^{1}}^{1}\,\extd \,q^{2}_{x}+ q_{x}^{1}\,\extd
\,q^{2}_{p^{1}_{x}}
     -q^{1}_{x}\,\extd \,p^{1}_{q^{2}_{x}})\wedge \extd p^{2}_{x}  &=&
      +\,x'(g;gn\so )\,\extd \,x'(g;gn\st )\wedge \extd \,x'(g;gn\sth )\cr
&&-\,
     x'(g;n\so g)\,\extd \,x'(g;n\st g)\wedge \extd \,x'(g;n\sth g)\ ,
\end{eqnarray*} where $n=[\![r_-,r_-]\!]$.  As $n$ is Adjoint
invariant, we have shown that
\begin{eqnarray*}
   J_{2}(x,x,x)
      & = & 4\,h_{p_{x}^{1}x}\wedge \extd p^{2}_{x}\,-\, 2\,p^{1}_{x}\,
\extd \,h_{p_{x}^{2}x}\,-\,2\,[\extd
x,h_{xx}]_\bullet/\hbar\,+\,[x,\extd \,h_{xx}]_\bullet/\hbar \,+\,
O(\hbar) \cr &=&
       4\,h_{p_{x}^{1}x}\wedge \extd p^{2}_{x}\,+\, 2\,p^{2}_{x}\,
\extd \,h_{p_{x}^{1}x}\,-\,2\,\extd
\,[x,h_{xx}]_\bullet/\hbar\,+\,3\,[x,\extd \,h_{xx}]_\bullet/\hbar
  \,+\,
O(\hbar) \ . \end{eqnarray*} Next using the equation
\[
D_{(g;\omega\st (g))}\,{\hat L}_{x}(v)\,-\, {\hat
L}_{p^{2}_{x}}(v)\,=\, x'(g;\omega\st (g)v-D_{(g;gv)}\omega\st
(g))\ ,
\]
we can write
\begin{eqnarray*}
h'_{xx}(g;\omega\st (g)) &=& 2\,\hat\Xi(D_{(g;\omega\st
(g))}\,{\hat L}_{x},{\hat L}_{x}) \cr &=& 2\, h_{p^{2}_{x}x}\,+\,
2\,\hat\Xi(v\mapsto x'(g;\omega\st (g)v-D_{(g;gv)}\omega\st
(g)),{\hat L}_{x})\ .\cr p^{1}_{y}(g)\,h'_{xx}(g;\omega\st (g))
&=& 2\,p^{1}_{y}(g) \, h_{p^{2}_{x}x}(g)\,+\, 2\, y'(g;g{r_-}\so
)\,\hat\Xi(v\mapsto x'(g;g[{r_-}\st ,v]),{\hat L}_{x}) \cr &=&
  2\,p^{1}_{y}(g)
\, h_{p^{2}_{x}x}(g)\,-\,2\,{\hat L}_{y}(g)({r_-}\so
)\,\hat\Xi(\ad^{*}_{{r_-}\st }\, {\hat L}_{x},{\hat L}_{x})\ .
\end{eqnarray*} Then, using the ad-invariance of $\hat\Xi$,
\begin{eqnarray*}
[y,h_{xx}]_\bullet/\hbar &=& x'(g;\omega\so
(g))\,h'_{xx}(g;\omega\st (g))\,+\, \Xi({\hat L}_{x},h_{xx})
\,+\,O(\hbar)  \cr
  &=& 2\,p^{1}_{y}(g)
\, h_{p^{2}_{x}x}(g)\,+\,\hat\Xi({\hat L}_{y},h_{xx}) \cr && +\,
{\hat L}_{y}(g)({r_-}\so )\big(\ad^{*}_{{r_-}\st } \hat\Xi({\hat
L}_{x},{\hat L}_{x})\,-\,2\, \hat\Xi(\ad^{*}_{{r_-}\st }{\hat
L}_{x},{\hat L}_{x})\big) \,+\,O(\hbar) \cr &=&2\,p^{1}_{y}(g) \,
h_{p^{2}_{x}x}(g)\,+\,\hat\Xi({\hat L}_{y},h_{xx})\,+\,O(\hbar) \
. \end{eqnarray*} This gives
\begin{eqnarray*}
   J_{2}(x,x,x)
      &=&
      4\,h_{p_{x}^{1}x}\wedge \extd p^{2}_{x}\,+\,
      2\,p^{2}_{x}\, \extd \,h_{p_{x}^{1}x}\,-\,2\,\extd
\,\Big(2\,p^{1}_{x}(g) \,
      h_{p^{2}_{x}x}(g)  \cr && \,+\,\hat\Xi({\hat L}_{x},h_{xx})
      \Big)\,+\,3\,[x,\extd \,h_{xx}]_\bullet/\hbar\,+\,O(\hbar) \cr
      &=&
      4\,h_{p_{x}^{1}x}\wedge \extd p^{2}_{x}\,+\,
      2\,p^{2}_{x}\, \extd \,h_{p_{x}^{1}x}\,+\,4\,\extd \,(p^{2}_{x}(g) \,
      h_{p^{1}_{x}x}(g)) \cr &&-\,2\,\extd \,\hat\Xi({\hat L}_{x},h_{xx})
      \,+\,3\,[x,\extd \,h_{xx}]_\bullet/\hbar\,+\,O(\hbar) \cr
      &=&
      6\,p^{2}_{x}\, \extd \,h_{p_{x}^{1}x}
      \,-\,2\,\extd \,\hat\Xi({\hat L}_{x},h_{xx})
      \,+\,3\,[x,\extd \,h_{xx}]_\bullet/\hbar\,+\,O(\hbar)\ .
\end{eqnarray*} Next we find, where $E={\hat
L}_{x}(g)({r_-}\so )\,\hat\Xi(\ad^{*}_{{r_-}\st }\, {\hat
L}_{x},{\hat L}_{x})$,
\begin{eqnarray*}
&&\nquad\quad p^{1}_{x}(g)\,h'_{xx}(g;\omega\st (g))
=2\,p^{1}_{x}(g) \, h_{p^{2}_{x}x}(g)\,-\,2\,E\ ,\cr &&\nquad\quad
p_{x}^{1}\,\extd \,D_{(g;\omega\st (g))}\,h_{xx} = \extd (\,
p_{x}^{1}\,D_{(g;\omega\st (g))}\,h_{xx}\,)\,-\,\extd
\,p_{x}^{1}\wedge D_{(g;\omega\st (g))}\,h_{xx} \cr &=& 2\, \extd
(\, p^{1}_{x}(g) \, h_{p^{2}_{x}x}(g)\,)\,-\,2\,\extd
\,E\,-\,\extd \,p_{x}^{1}\wedge D_{(g;\omega\st (g))}\,h_{xx} \cr
&=& \extd \, p^{1}_{x}(g) \wedge(2\,
h_{p^{2}_{x}x}(g)-D_{(g;\omega\st (g))}\,h_{xx}) \,+\, 2\,
p^{1}_{x}(g) \,\extd \, h_{p^{2}_{x}x}(g)\,-\,2\,\extd \,E \cr
  &=&
2\, \extd \,({\hat L}_{x}(g)({r_-}\so ))\wedge
  \hat\Xi(\ad^{*}_{{r_-}\st }\,
{\hat L}_{x},{\hat L}_{x}) \,+\, 2\, p^{1}_{x}(g) \,\extd \,
h_{p^{2}_{x}x}(g)\,-\,2\,\extd \,E \cr &=&
  2\, p^{1}_{x}(g) \,\extd \, h_{p^{2}_{x}x}(g)\,-\,2\,{\hat
L}_{x}(g)({r_-}\so )\,
  \extd \, \hat\Xi(\ad^{*}_{{r_-}\st }\,
{\hat L}_{x},{\hat L}_{x})\ . \end{eqnarray*} Using this with
\ref{rope} and the ad-invariance of $\hat\Xi$ again,
\begin{eqnarray*}
[x,\extd h_{xx}]_\bullet(v,w)/\hbar &=& 2\, p^{1}_{x}(g) \,\extd
\, h_{p^{2}_{x}x}(g)(v,w) \,+\, (\hat L_{x}\tens \extd
\,h_{xx})(\hat\Xi(v)\tens w-\hat\Xi(w)\tens v) \cr && \kern-20pt
+\,\hat L_{x}({r_-}\so )\,\extd (\ad^{*}_{{r_-}\st }(h_{xx})-2\,
\hat\Xi(\ad^{*}_{{r_-}\st }\, \hat L_{x},\hat L_{x}))(g)(v,w)
\,+\,O(\hbar) \cr &=& 2\, p^{1}_{x}(g) \,\extd \,
h_{p^{2}_{x}x}(g)(v,w)\,+\,G(v,w) \,+\,O(\hbar)\ , \end{eqnarray*}
where $G(v,w)=(\hat L_{x}\tens \extd \,h_{xx})(\hat\Xi(v)\tens
w-\hat\Xi(w)\tens v)$.  This results in
\begin{eqnarray*}
     J_{2}(x,x,x) &=& 3\,G
     \,-\,2\, \extd \,\hat\Xi(\hat L_{x},h_{xx})\ .
\end{eqnarray*} Remember that
$h_{xx}(g)(v)=x'(g;g\hat\Xi\so (v))\,x'(g;g\hat\Xi\st (v))$, so
using Lemma~\ref{string},
\begin{eqnarray*}
     \extd \,h_{xx}(g)(w,v) &=&
     2\,x''(g;g\hat\Xi\so (v),gw)\,x'(g;g\hat\Xi\st (v))\,+\,2\,
     x'(g;gw\,\hat\Xi\so (v))\,x'(g;g\hat\Xi\st (v))\cr &&-\,
     2\,x''(g;g\hat\Xi\so (w),gv)\,x'(g;g\hat\Xi\st (w))\,-\,2\,
     x'(g;gv\,\hat\Xi\so (w))\,x'(g;g\hat\Xi\st (w))\cr &&-\,
     x'(g;g\hat\Xi\so ([w,v]))\,x'(g;g\hat\Xi\st ([w,v]))\ .
\end{eqnarray*} Then, using the symmetry of $\hat\Xi$ and an
independent copy $\hat\Xi'$,
\begin{eqnarray*}
G(v,w) &=& \hat L_x(\hat\Xi\so (v))\, \extd \,h_{xx}(\hat\Xi\st
(v),w)\,-\,
  \hat L_x(\hat\Xi\so (w))\, \extd \,h_{xx}(\hat\Xi\st (w),v) \cr
&=& 2\,x''(g;g\hat\Xi'\so (\hat\Xi\st (w)),gv)\,{\hat L}_{x}
(\hat\Xi'\st (\hat\Xi\st (w)))\,{\hat L}_{x}(\hat\Xi\so (w))\cr
&&-\, 2\,x''(g;g\hat\Xi'\so (\hat\Xi\st (v)),gw)\,{\hat L}_{x}
(\hat\Xi'\st (\hat\Xi\st (v)))\,{\hat L}_{x}(\hat\Xi\so (v))\cr
&&+\, 2\, {\hat L}_{x}(\hat\Xi\st (v)\,\hat\Xi'\so (w))\, {\hat
L}_{x}(\hat\Xi'\st (w))\, {\hat L}_{x}(\hat\Xi\so (v))\\
&& -\, 2\, {\hat L}_{x}(w\,\hat\Xi'\so (\hat\Xi\st (v)))\, {\hat
L}_{x}(\hat\Xi'\st (\hat\Xi\st (v)))\, {\hat L}_{x}(\hat\Xi\so
(v)) \cr &&-\, 2\, {\hat L}_{x}(\hat\Xi\st (w)\,\hat\Xi'\so (v))\,
{\hat
L}_{x}(\hat\Xi'\st (v))\, {\hat L}_{x}(\hat\Xi\so (w))\\
&& +\, 2\, {\hat L}_{x}(v\,\hat\Xi'\so (\hat\Xi\st (w)))\, {\hat
L}_{x}(\hat\Xi'\st (\hat\Xi\st (w)))\, {\hat L}_{x}(\hat\Xi\so
(w)) \cr && -\, {\hat L}_{x}(\hat\Xi'\so ([\hat\Xi\st (v),w])) \,
{\hat L}_{x}(\hat\Xi'\st ([\hat\Xi\st (v),w])) \, {\hat
L}_{x}(\hat\Xi\so (v)) \cr  && +\, {\hat L}_{x}(\hat\Xi'\so
([\hat\Xi\st (w),v])) \, {\hat L}_{x}(\hat\Xi'\st ([\hat\Xi\st
(w),v])) \, {\hat L}_{x}(\hat\Xi\so (w)) \cr &=&
2\,D_{(g;gv)}\Big({\hat L}_{x}(\hat\Xi'\so (\hat\Xi\st
(w)))\Big)\, {\hat L}_{x} (\hat\Xi'\st (\hat\Xi\st (w)))\,{\hat
L}_{x}(\hat\Xi\so (w))\cr &&-\, 2\,D_{(g;gw)}\Big({\hat
L}_{x}(\hat\Xi'\so (\hat\Xi\st (v)))\Big)\, {\hat L}_{x}
(\hat\Xi'\st (\hat\Xi\st (v)))\,{\hat L}_{x}(\hat\Xi\so (v)) \cr
&&+\, 2\, {\hat L}_{x}([\hat\Xi\st (v),\hat\Xi'\so (w)])\, {\hat
L}_{x}(\hat\Xi'\st (w))\, {\hat L}_{x}(\hat\Xi\so (v))\cr && -\,
{\hat L}_{x}(\hat\Xi'\so ([\hat\Xi\st (v),w])) \, {\hat
L}_{x}(\hat\Xi'\st ([\hat\Xi\st (v),w])) \, {\hat
L}_{x}(\hat\Xi\so (v)) \cr  && +\, {\hat L}_{x}(\hat\Xi'\so
([\hat\Xi\st (w),v])) \, {\hat L}_{x}(\hat\Xi'\st ([\hat\Xi\st
(w),v])) \, {\hat L}_{x}(\hat\Xi\so (w))\ .
\end{eqnarray*} Using Lemma~\ref{string},
\begin{eqnarray*}
\extd \,\hat\Xi(\hat L_{x},h_{xx})(v,w) &=& D_{(g;gv)}(\hat
L_{x}(\hat\Xi\so (w))\,h_{xx}(\hat\Xi\st (w)))\\
&& \kern-20pt -\, D_{(g;gw)}(\hat L_{x}(\hat\Xi\so
(v))\,h_{xx}(\hat\Xi\st (v))) -\, \hat L_{x}(\hat\Xi\so
([v,w]))\,h_{xx}(\hat\Xi\st ([v,w]))\ ,
\end{eqnarray*} so
\begin{eqnarray*}
      {\hat L}_{x}(\hat\Xi\so (w))\,h_{xx}(\hat\Xi\st (w)) &=&
      {\hat L}_{x}(\hat\Xi'\so \hat\Xi\st (w))\,
      \,{\hat L}_{x}(\hat\Xi'\st \hat\Xi\st (w))\,{\hat L}_{x}(\hat\Xi\so
(w))\ ,\cr
      D_{(g;gv)}({\hat L}_{x}(\hat\Xi\so (w))\,h_{xx}(\hat\Xi\st (w)))&=&
      2\,D_{(g;gv)}({\hat L}_{x}(\hat\Xi'\so \hat\Xi\st (w)))\,
      \,{\hat L}_{x}(\hat\Xi'\st \hat\Xi\st (w))\,{\hat L}_{x}(\hat\Xi\so
(w))\cr &&+\,
      {\hat L}_{x}(\hat\Xi'\so \hat\Xi\st (w))\,
      \,{\hat L}_{x}(\hat\Xi'\st \hat\Xi\st (w))\,D_{(g;gv)}({\hat
L}_{x}(\hat\Xi\so (w))) \ .
\end{eqnarray*} Now we get
\begin{eqnarray*}
     J_{2}(x,x,x)(v,w) &=&
     2\,D_{(g;gv)}\Big({\hat L}_{x}(\hat\Xi'\so (\hat\Xi\st (w)))\Big)\,
{\hat L}_{x} (\hat\Xi'\st (\hat\Xi\st (w)))\,{\hat
L}_{x}(\hat\Xi\so (w))\cr &&-\, 2\,D_{(g;gw)}\Big({\hat
L}_{x}(\hat\Xi'\so (\hat\Xi\st (v)))\Big)\, {\hat L}_{x}
(\hat\Xi'\st (\hat\Xi\st (v)))\,{\hat L}_{x}(\hat\Xi\so (v)) \cr
&&+\, 6\, {\hat L}_{x}([\hat\Xi\st (v),\hat\Xi'\so (w)])\, {\hat
L}_{x}(\hat\Xi'\st (w))\, {\hat L}_{x}(\hat\Xi\so (v))\cr && -\,
3\,{\hat L}_{x}(\hat\Xi'\so ([\hat\Xi\st (v),w])) \, {\hat
L}_{x}(\hat\Xi'\st ([\hat\Xi\st (v),w])) \, {\hat
L}_{x}(\hat\Xi\so (v)) \cr  && +\,3\, {\hat L}_{x}(\hat\Xi'\so
([\hat\Xi\st (w),v])) \, {\hat L}_{x}(\hat\Xi'\st ([\hat\Xi\st
(w),v])) \, {\hat L}_{x}(\hat\Xi\so (w)) \cr &&-\,2\,{\hat
L}_{x}(\hat\Xi'\so \hat\Xi\st (w))\,
      \,{\hat L}_{x}(\hat\Xi'\st \hat\Xi\st (w))\,D_{(g;gv)}({\hat
L}_{x}(\hat\Xi\so (w))) \cr &&+\,2\,{\hat L}_{x}(\hat\Xi'\so
\hat\Xi\st (v))\, \,{\hat L}_{x}(\hat\Xi'\st \hat\Xi\st
(v))\,D_{(g;gw)}({\hat L}_{x}(\hat\Xi\so (v)))\cr &&+\,2\,{\hat
L}_{x}(\hat\Xi\so ([v,w]))\,h_{xx}(\hat\Xi\st ([v,w]))\ .
\end{eqnarray*} By ad-invariance of $\hat\Xi'$,
\begin{eqnarray*}
\hat\Xi'([w,\hat\Xi\st (v)]) &=& [\hat\Xi'\so (w),\hat\Xi\st (v)]
\tens \hat\Xi'\st (w)\,+\,\hat\Xi'\so (w)\tens [\hat\Xi'\st
(w),\hat\Xi\st (v)]\ ,\cr {\hat L}_{x}(\hat\Xi\so
([v,w]))\,h_{xx}(\hat\Xi\st ([v,w])) &=& {\hat L}_{x}([\hat\Xi\so
(v),w])\, {\hat L}_{x}(\hat\Xi'\so \hat\Xi\st (v)) \, {\hat
L}_{x}(\hat\Xi'\st \hat\Xi\st (v)) \cr && +\, 2\, {\hat
L}_{x}(\hat\Xi\so (v))\, {\hat L}_{x}([\hat\Xi\st (v),\hat\Xi'\so
(w)])\, {\hat L}_{x}(\hat\Xi'\st (w))\ , \end{eqnarray*} so we
rewrite
\begin{eqnarray*}
     J_{2}(x,x,x)(v,w) &=&
     2\,D_{(g;gv)}\Big({\hat L}_{x}(\hat\Xi'\so (\hat\Xi\st (w)))\Big)\,
{\hat L}_{x}
     (\hat\Xi'\st (\hat\Xi\st (w)))\,{\hat L}_{x}(\hat\Xi\so (w))\cr &&-\,
     2\,D_{(g;gw)}\Big({\hat L}_{x}(\hat\Xi'\so (\hat\Xi\st (v)))\Big)\,
{\hat L}_{x}
     (\hat\Xi'\st (\hat\Xi\st (v)))\,{\hat L}_{x}(\hat\Xi\so (v)) \cr &&-\,
2\,
     {\hat L}_{x}([\hat\Xi\st (v),\hat\Xi'\so (w)])\, {\hat
L}_{x}(\hat\Xi'\st (w))\,
     {\hat L}_{x}(\hat\Xi\so (v))\cr
     &&-\,2\,{\hat L}_{x}(\hat\Xi'\so \hat\Xi\st (w))\,
     \,{\hat L}_{x}(\hat\Xi'\st \hat\Xi\st (w))\,D_{(g;gv)}({\hat
L}_{x}(\hat\Xi\so (w)))
     \cr &&+\,2\,{\hat L}_{x}(\hat\Xi'\so \hat\Xi\st (v))\,
     \,{\hat L}_{x}(\hat\Xi'\st \hat\Xi\st (v))\,D_{(g;gw)}({\hat
L}_{x}(\hat\Xi\so (v)))\cr
     &&+\,{\hat L}_{x}([\hat\Xi\so (v),w])\, {\hat L}_{x}(\hat\Xi'\so
\hat\Xi\st (v)) \, {\hat L}_{x}(\hat\Xi'\st \hat\Xi\st (v)) \cr
&&-\, {\hat L}_{x}([\hat\Xi\so (w),v])\, {\hat L}_{x}(\hat\Xi'\so
\hat\Xi\st (w)) \, {\hat L}_{x}(\hat\Xi'\st \hat\Xi\st (w))
\end{eqnarray*} Setting $e(v)=\hat\Xi\so (v)\tens
\hat\Xi'\so \hat\Xi\st (v)\tens
    \hat\Xi'\st \hat\Xi\st (v)\,-\, \hat\Xi'\so \hat\Xi\st (v)\tens
    \hat\Xi\so (v)\tens \hat\Xi'\st \hat\Xi\st (v)$,
\begin{eqnarray*}
     J_{2}(x,x,x)(v,w) &=&
     2\,D_{(g;gw)}\big({\hat L}_{x}(e\so (v))\big)\, L_x(e\st (v))\,
L_x(e\sth (v))  \cr &&  -\, 2\,D_{(g;gv)}\big({\hat L}_{x}(e\so
(w))\big)\, L_x(e\st (w))\, L_x(e\sth (w))
  \cr &&-\, 2\,
     {\hat L}_{x}([\hat\Xi\st (v),\hat\Xi'\so (w)])\, {\hat
L}_{x}(\hat\Xi'\st (w))\,
     {\hat L}_{x}(\hat\Xi\so (v))\cr
     && +\,\hat L_x([e\so (v),w])\,\hat L_x(e\st (v))\,\hat L_x(e\sth (v))
     \\ &&-\, \hat L_x([e\so (w),v])\,\hat L_x(e\st (w))\,\hat L_x(e\sth (w))
\cr
  &&+\,{\hat L}_{x}([\hat\Xi'\so \hat\Xi\st (v),w])\, {\hat L}_{x}(\hat\Xi\so
(v)) \, {\hat L}_{x}(\hat\Xi'\st \hat\Xi\st (v)) \cr &&-\, {\hat
L}_{x}([\hat\Xi'\so \hat\Xi\st (w),v])\, {\hat L}_{x}(\hat\Xi\so
(w)) \, {\hat L}_{x}(\hat\Xi'\st \hat\Xi\st (w)) \ .
\end{eqnarray*} Now, using symmetry and ad-invariance,
\begin{eqnarray*}
&& 2\,{\hat L}_{x}([\hat\Xi'\so \hat\Xi\st (v),w]) \, {\hat
L}_{x}(\hat\Xi'\st \hat\Xi\st (v)) \\
&&\quad= {\hat L}_{x}([\hat\Xi'\so \hat\Xi\st (v),w]) \, {\hat
L}_{x}(\hat\Xi'\st \hat\Xi\st (v))\,+\, {\hat L}_{x}(\hat\Xi'\so
\hat\Xi\st (v))\, {\hat
L}_{x}([\hat\Xi'\st \hat\Xi\st (v),w])\\
&&\quad={\hat L}_{x}(p\so)\,{\hat L}_{x}(p\st)\ , \end{eqnarray*}
where
\begin{eqnarray*}
p\so\tens p\st &=&
\hat\Xi'([\hat\Xi\st(v),w])\,=\,[\hat\Xi\st(v),\hat\Xi'\so(w)]\tens
\hat\Xi'\st(w) \,+\,\hat\Xi'\so(w)\tens
[\hat\Xi\st(v),\hat\Xi'\st(w)]\ . \end{eqnarray*} It follows that
\begin{eqnarray*}
{\hat L}_{x}([\hat\Xi'\so \hat\Xi\st (v),w]) \, {\hat
L}_{x}(\hat\Xi'\st \hat\Xi\st (v)) &=& {\hat L}_{x}([\hat\Xi\st
(v),\hat\Xi'\so (w)]) \, {\hat L}_{x}(\hat\Xi'\st (w))\ ,
\end{eqnarray*} and so
\begin{eqnarray*}
     J_{2}(x,x,x)(v,w) &=&
     2\,D_{(g;gw)}\big({\hat L}_{x}(e\so (v))\big)\, \hat L_x(e\st (v))\,
     \hat L_x(e\sth (v))  \cr &&  -\,
2\,D_{(g;gv)}\big({\hat L}_{x}(e\so (w))\big)\, \hat L_x(e\st
(w))\, \hat L_x(e\sth (w))
  \cr
     && +\,\hat L_x([e\so (v),w])\,\hat L_x(e\st (v))\,\hat L_x(e\sth (v))
     \\
     &&-\, \hat L_x([e\so (w),v])\,\hat L_x(e\st (w))\,\hat L_x(e\sth (w))\ .
\end{eqnarray*} Using the first super Jacobi result $
    e(v)  \,=\, -\,n\so \tens
    n\st \tens [v,n\sth]$ as stated. \eproof

From Section~2 we know that if $J_1,J_2$ vanish then the third
super Jacobi identity also holds. Putting several of the above
results together we find as a special case:

\begin{theorem}\label{cantilde} Every quastriangular Poisson-Lie group
has a compatible bicovariant preconnection given by
$\hat\Xi=0$. The corresponding $\tilde\gamma$ is given by \[
\tilde\gamma(a,s)(g)=\{a,s\}(g)+a'(g;gr_-\so)\ad^*_{r_-\st}(s(g))\]
for $a\in C^\infty(G)$ and $s\in C^\infty(G,\cg^*)$. If $\cg$ is
semisimple then the curvature vanishes if and only if the Lie
bialgebra is triangular, and in this case all super Jacobi
identities $J_i=0$ hold.
\end{theorem}
\proof Clearly $\hat\Xi=0$ is symmetric and $\ad^*$-invariant,
hence by Corollary~\ref{uuttii} and Proposition~\ref{hatXiad}
defines a bicovariant preconnection. To work out what it
looks like we have only to work backwards from the definitions of
$\hat\Xi,\Xi$ and Definition~\ref{connen} to find the result
stated. If $\hat\Xi\ne 0$ we would need an additional term
$\hat\Xi(\hat L_x (g),s(g))$. For the second part, the Lie
bialgebra is triangular precisely when
$n=[[r_-,r_-]]=[r_{+12},r_{+23}]$ vanishes, which is if and only
if $r_+=0$ since this is either zero or nondegenerate in the
semisimple case. Also in this case $e(v)=0$ for all $v$ if and
only if $n=0$, which is the case for the curvature obstruction to
$J_1$ to vanish according to Proposition \ref{qbicj1} and the
interpretation in Section~2. In this case $J_2$ also vanishes by
the preceding proposition and hence $J_3$ by
Proposition~\ref{jac3}. \eproof

\begin{propos}\label{cangamma}
The above canonical choice of compatible preconnection
$\gamma$ on a quasitriangular Poisson-Lie group is given
explicitly by \[ \gamma(a, \tau)=(a\ra r_-\so)(\tau\ra
r_-\st)-(r_-\so\la a)(r_-\st\la \tau)\] where $\ra$ is the right
action on on functions or 1-forms corresponding to $\Delta,
\Delta_R$ respectively. Its curvature and torsion are
\begin{eqnarray*} R(\hat{x},\hat{y})\tau&=&(n\so\la x)(n\st\la
y)(n\sth\la
\tau)-(x\ra n\so)(y\ra n\st)(\tau\ra n\sth)\\
  T(\hat x,\hat y)(\extd z)&=&(m\so\la x)(m\st\la y)(m\sth\la z)-(x\ra m\so)(y\ra m\st)(z\ra
m\sth)\end{eqnarray*} where $n=[[r_-,r_-]]$ and
$m=[r_{-13},r_{-23}]$.
\end{propos} \proof Here the right actions are defined in the same way as we
did for $\la$ at the start of Section~4.1, so $(a\ra
g)(h)=a(hg^{-1})$ on functions and $(\tau\ra
g)(h)=R^*_{g^{-1}}(\tau(hg^{-1}))$ on forms, and we use the
infinitesimal versions. Unwinding our definition of
$\tilde\gamma$, we have
\begin{eqnarray*}
\gamma(a,\tau)(g)&=&L^*_{g^{-1}}(\tilde\gamma(a,\tilde\tau)(g))\\
&=&L^*_{g^{-1}}\left(a'(g;gr_-\so)(\tilde\tau'(g;gr_-\st)
+\ad^*_{r_-\st}(\tilde\tau(g)))-a'(g;r_-\so g
)\tilde\tau'(g;r_-\st g)\right)\\ &=&a'(g;r_-\so g)(r_-\st\la
\tau)(g)-a'(g;gr_-\so)(\tau\ra r_-\st)(g)
\end{eqnarray*}
as stated. Here \[ a'(g;gv)=(L_{*v}a)(g)=-(a\ra v)(g),\quad
a'(g;vg)=(R_{*v}a)(g)=-(v\la a)(g)\] for $v\in\cg$. On forms, we
similarly have
\[ L_{g^{-1}}^*(\tilde\tau'(g;vg))
=L^*_{g^{-1}}{\extd \over \extd t}|_0\tilde\tau(e^{tv}g)={\extd
\over \extd t}|_0L^*_{e^{tv}}(\tau(e^{tv}g))={\extd \over \extd
t}|_0(e^{-tv}\la \tau)(g)=-(v\la\tau)(g)\]
\begin{eqnarray*}L_{g^{-1}}^*(\tilde\tau'(g;gv)
+\ad^*_v(\tilde\tau(g)))&=&L^*_{g^{-1}} {\extd \over \extd
t}|_0\Ad^*_{e^{tv}}(\tilde\tau(ge^{tv}))={\extd \over \extd
t}|_0R^*_{e^{tv}}(\tau(ge^{tv}))\\&=&{\extd \over \extd
t}|_0(\tau\ra e^{-tv})(g)=-(\tau\ra v)(g).\end{eqnarray*} This is
equivalent to the computation in Lemma~\ref{upup} for exact forms.
Next, it is a useful check to compute the curvature directly form
the definition (\ref{partialR}), with result in line with
Proposition~\ref{qbicj1}. We use the canonical $\gamma$ and its
connection property (\ref{eqq1}) to write
\begin{eqnarray*}
R(\hat{x},\hat{y})\tau&=&\gamma(x,\gamma(y,\tau))
-\gamma(y,\gamma(x,\tau))-\gamma(\{x,y\},\tau)\\
&=&\gamma(x,(y\ra r_-\so)(\tau\ra r_-\st)-(r_-\so\la y)(r_-\st\la
\tau))\\
&&-(x\leftrightarrow y)-(\{x,y\}\ra r_-\so)(\tau\ra
r_-\st)+(r_-\so\la\{x,y\})(r_-\st\la\tau).\end{eqnarray*} We then
insert the formulae for $\gamma$ and $\{\ ,\ \}$ (as discussed
below) and expand $\la,\ra$ as derivations to obtain 24 terms.
Cancelling 12 and using that $\la,\ra$ are mutually commuting
actions of $\cg$, and the antisymmetry of $r_-$ then gives the
result. Finally, using the expression for torsion as in
Proposition~\ref{tor} we have
\begin{eqnarray*} \<\extd z, T(\hat x,\hat y)\>&=&\<\hat
x,\gamma(y,\extd z)\>-\<\hat y,\gamma(x,\extd z)\>\\
&=&\<\hat x,(y\ra r_-\so) \extd (z\ra r_-\st)-(r_-\so\la y)\extd
(r_-\st\la z)\>-(x\leftrightarrow y)\\
&=&(y\ra r_-\so)\{x,z\ra r_-\st\}-(r_-\so\la y)\{x,r_-\st\la
z\}-(x\leftrightarrow y)\end{eqnarray*} where we used $\gamma$ on
exact forms as discussed below. We then expand out the Poisson
bracket and cancel terms. One may similarly compute $T(\hat x,\hat
y)(\tau)$ in general. \eproof

Since $\extd$ commutes with the actions, or by Lemma~\ref{upup},
we have in particular \eqn{canexact}{\gamma(a,\extd b)=( a\ra
r_-\so) \extd (b\ra r_-\st)-(r_-\so\la a)\extd (r_-\st\la b)}
while in the same notation, \eqn{poirmin}{ \{a,b\}=( a\ra r_-\so)
(b\ra r_-\st)-(r_-\so\la a) (r_-\st\la b)} coincides with the
usual formula on a quasitriangular Poisson-Lie group (because one
could equally well put $r$ here since $r_+$ is $\ad$-invariant).
Then antisymmetry of $r_-$ and the Leibniz rule confirms that the
canonical connection is compatible in the sense of (\ref{eqq2}).
Also, the torsion is consistent with Proposition~\ref{tor} since
on cyclic sum, $[r_{-13},r_{-23}]$ becomes replaced by
$[[r_-,r_-]]$ and this now gives zero as the statement that
(\ref{poirmin}) obeys the Jacobi identity (usually this is done
via $[[r,r]]=0$). By contrast, the action of $[[r_-,r_-]]$ in the
curvature, even on exact forms, is not trivial. We note also that
the same formula as for the canonical connection but with $r$ in
place of $r_-$ is what one obtains by semiclassicalising the
`quantum Lie functor' construction in \cite{GomMa} for a canonical
bicovariant quantum differential calculus on any coquasitriangular
Hopf algebra. There the coquasitriangular structure $\mathcal R$
or `universal R-matrix' plays the role of $r$. Unfortunately, that
functor only gives nontrivial answers in the triangular case and
now we can see why: only in this case is $r=r_-$ so that it
coincides with our above result. Moreover, only this case
corresponds to $J_1=0$ and hence to a bimodule on quantization.

Finally, let us use Corollary~4.7 and Proposition~4.16 to study
the entire moduli space of bicovariant semiclassical calculi.
These results tell us that the moduli space is an affine space
with the above canonical preconnection as reference and
others given by the vector space of $\Ad$-invariant symmetric maps
$\hat\Xi:\cg\to {\rm Sym}^2(\cg)$. Such maps may be classified for
$\cg$ reductive using Kostant's results on harmonic functions in
\cite{Kos}. In particular:

\begin{theorem}\label{moduli} If $\cg$ is simple and not $sl_n$, $n>2$, the
moduli space is a point, i.e. the canonical preconnection
$\hat\Xi=0$ is the unique bicovariant one. If $\cg=sl_n$, $n>2$,
the moduli space is a 1-parameter family given by $\hat\Xi$ a
multiple of the split cubic Casimir in $\cg\tens \cg\tens \cg$.
\end{theorem} \proof Since $\cg$ is simple $\hat \Xi$ is zero or an
inclusion, so we need the multiplicity of $\cg$ in the symmetric
tensor square ${\rm Sym}^2(\cg)$. We freely use the Ad-invariant
Killing form to identify $\cg$ and $\cg^*$ for our purposes, so
this is equivalent to the multiplicity of $\cg$ in the symmetric
polynomials $S^2(\cg)$ of degree 2 on $\cg$. From \cite{Kos},
$S(\cg)=J\tens {\rm Harm}$ as a vector space, where $J$ denotes
the invariant polynomials and ${\rm Harm}$ the harmonic ones.
Hence $S^2(\cg)=J^2\oplus (J^1\tens {\rm Harm}^1)\oplus {\rm
Harm}^2$. $J$ is well-known to be generated by functions $\{u_i\}$
of degrees $m_i+1$ where $m_i$ are the `exponents' of the Lie
algebra. Clearly $J^2$ is 1-dimensional (spanned by the Killing
form) and $J^1=0$, so we need only to classify embeddings of $\cg$
in ${\rm Harm}^2$. From \cite{Kos} (see Theorem~0.11) the
component of ${\rm Harm}$ transforming as a given highest weight
representation $\lambda$ consists of $l_\lambda$ copies with
certain degrees $m_i(\lambda)$. In the case of the adjoint
representation these reduce to the usual rank $l$ and exponents
$m_i$.  Finally, we turn to tables and find that only $sl_n$,
$n>2$ has an exponent with value $2$ (i.e. a cubic Casimir). Other
simple Lie algebras have nothing in {\rm Harm} of degree 2 and
hence have $\hat \Xi=0$. In the case of $sl_n$, $n>2$ there is one
exponent with value 2  (namely $m_2$) so there is a single copy of
$\cg$ in the decomposition of $S^2(\cg)$. It is necessarily given
by the corresponding unique totally symmetric invariant element
$u_2\in J^3$, viewed as a map $\cg\to S^2(\cg)$. From another
point of view this corresponds to the cubic Casimir in the
enveloping algebra, or the unique totally symmetric `cubic split
Casimir' in ${\rm Sym}^3(\cg)$. Thus up to normalisation, $\hat
\Xi$ is necessarily this element viewed as $\cg^*\to {\rm
Sym}^2(\cg)$ by evaluation in one input and converted to $\cg\to
{\rm Sym}^2(\cg)$ with the help of the Killing form. \eproof

We see that for all but the $sl_n$ series, the corresponding
Poisson-Lie group with any fixed strictly quasitriangular
structure (such as the standard one) admits a unique semiclassical
bicovariant calculus, and its preconnection has nonzero
curvature. Therefore it cannot be quantised even at a bimodule
level, i.e. there can be no first order bicovariant different
calculus on the standard quantum groups  of the same dimension as
the classical one. This agrees with what is known from quantum
group theory by other means\cite{Ma:cla}. The 1-parameter family
for $sl_n$, $n>2$ at least generically also has curvature from
Proposition~4.16, and hence the same problem.

\subsection{Example: canonical connection on $SU_2$}

We take the same Poisson-Lie group as in Section~4.3 and the
notations there. This time there is no known quantum calculus to
semiclassicalise but rather we use our results above. In this case
one can see directly that $sl_2$ is not to be found in the
symmetric part of $sl_2\tens sl_2$ (which is the $1\oplus 5$
dimensional representation under $\ad$), hence  the only
bicovariant preconnection by Proposition \ref{hatXiad} is
$\hat\Xi=0$, as per Theorem~\ref{moduli}.

In this case Theorem~\ref{cantilde} says that $
\tilde\gamma(x,\tilde\tau^i)=(L_{*r_-\so}x)\ad^*_{r_-\st}\tilde\tau^i$
since the $\tilde\tau^i$ are constant. Also for this reason, we
have the same formula without all the tildes, with the $\tau^i$
transforming in the same way under $\ad^*$ as the $\tilde\tau^i$.
Putting in the form of $r_-$ in our case, we have
\[ \gamma(x,\tau^i)={1\over
2}((\del_+x) \, \ad^*_{e_-}(\tau^i)-(\del_- x)\,
\ad^*_{e_+}(\tau^i)).\] Computing $\ad^*_{e_i}$ we find
\[ \gamma(x,\tau^+)=-(\del_- x)\tau^3,\quad
\gamma(x,\tau^-)=-(\del_+x)\tau^3,\quad \gamma(x,\tau^3)={1\over
2}((\del_+x)\tau^++(\del_-x)\tau^-).\] As we expect from the
Proposition~\ref{cangamma}, this compatible preconnection has
both curvature and torsion. Hence by Sections~2 and~3 it cannot be
quantized to an honest bimodule first order quantum calculus let
alone higher forms.

\section{Quasiassociative exterior algebras}

Here we give a setting for the quantisation of Poisson-Lie groups
where the preconnection has nonzero curvature, such as the
canonical one in Theorem~\ref{cantilde} for the strictly
quasitriangular case. This case  includes the standard
Drinfeld-Sklyanin Poisson-Lie groups as demonstrated for $SU_2$ in
Section~4.5. In such a case of curvature we know that the
quantisation must be nonassociative. We show now that this can be
controlled nicely in the setting of coquasi-Hopf algebras, where
there is a Drinfeld associator $\Phi$. We first work out the
relevant noncommutative differential geometry first and then
semiclassicalise it. As a result, we succeed to quantise the
canonical connection for such quasitriangular Poisson Lie groups.

For general constructions we work over a field $k$, which can also
(with care) be replaced by commutative ring such as $\C[[\hbar]]$
for the formal deformation theory. Our starting point is a class
of natural examples obtained by twisting as follows. If $H$ is a
Hopf algebra and $F:H\tens H\to k$ is a 2-cochain in the sense of
convolution invertible in the form $F(a\o\tens b\o)F^{-1}(a\t\tens
b\t)=\eps(a)\eps(b)$ (and similarly on the other side) and obeying
$F(1\tens a)=\eps(a)$ then one may define a new object $H_F$ with
modified product
  \eqn{twistH}{ a\bullet b=F(a\o\tens b\o) a\t b\t
F^{-1}(a\thr\tens b\thr),\quad\forall a,b\in H.} This is not
necessarily associative but rather \eqn{qassoc}{a\bullet(b\bullet
  c) =\Phi (a\o\tens b\o\tens c\o)(a\t\bullet b\t)
  \bullet c\t\Phi^{-1}(a\thr\tens b\thr\tens
  c\thr)}
where
\[ \Phi(a\tens b\tens c)=F(b\o\tens c\o)F(a\o\tens b\t c\t)F^{-1}(a\t
b\thr\tens c\thr) F^{-1}(a\thr\tens b\fo)
\]
is the coboundary of $F$ in some kind of nonAbelian
cohomology\cite{Ma:book}. This makes $H_F$ into an example of a
coquasi-Hopf algebra. By definition the latter is defined to be a
coalgebra which has a product obeying (\ref{qassoc}) with respect
to some 3-cocycle $\Phi$ (not necessarily of the coboundary form).
There should also be an antipode $S$ obeying certain axioms. The
axioms in a non-dual form are due to Drinfeld\cite{Dri:alm} who
showed what would in our setting be the following assertion: the
standard coquasitriangular quantum groups $C_\hbar(G)$ are (non
trivially) isomorphic to twists of the classical functions $C(G)$
(in some form) by a formal cochain $F=F_\hbar$. Although the $F$
needed here is not a cocycle so that $\Phi$ is not trivial, $\Phi$
turns out in this example to be cocentral and hence disappears
from (\ref{qassoc}) with the result that $C_\hbar(G)$ is an
ordinary Hopf algebra. We do not need the fully precise
formulation but only the idea and formulae to lowest order, which
we will provide.

The second idea is the Majid-Oeckl theorem for twisting of
bicovariant differential calculi on Hopf algebras. It is known
\cite{Brz:rem} that in the bicovariant case the  Woronowicz
exterior algebra $\Omega(H)$ as a whole forms a super-Hopf
algebra. It is just $H$ in degree 0 and in higher degree the
coproduct is generated by \[ \Delta=\Delta_L+\Delta_R\] in degree
1. We will use the notation
\eqn{bicomod}{(\Delta_L\tens\id)\Delta_R\tau=(\id\tens\Delta_R)
\Delta_L\tau=\tau\o \tens \tau\t\tens \tau\thr\in H\tens
\Omega^1\tens H} for any 1-form $\tau$. The theorem of
\cite{MaOec:twi} states that if $F$ is a 2-cocycle (so that $H_F$
is again a Hopf algebra) then $H_F$ has a natural bicovariant
differential calculus given by $\Omega(H_F)=\Omega(H)_F$, where we
twist $\Omega(H)$ as a super-Hopf algebra with $F$ extended by
zero on higher degrees. The twisted module and wedge product are
\begin{eqnarray}a\bullet \tau &=& F(a\o\tens \tau\o) a\t \tau\t
F^{-1}(a\thr\tens \tau\thr)\nonumber\\
\tau\bullet a&=&F(\tau\o\tens a\o) \tau\t a\t F^{-1}(\tau\thr\tens
a\thr)\nonumber \\ \label{twicalc} \tau\wedge_\bullet \eta&=&
F(\tau\o\tens \eta\o) \tau\t\wedge \eta\t F^{-1}(\eta\thr\tens
\eta\thr)\end{eqnarray} for $a\in H$ and 1-forms $\tau,\eta$. Note
that $(\Delta\tens\id)\Delta\tau$ contains three terms, but only
the middle one (\ref{bicomod}) enters in the above expressions
since $F$ pairs only in degree zero. The same applies in formulae
below.

We can ask what kind of object does one get if $F$ is not a
cocycle. According to the above, then $\Omega(H_F)$ given by
twisting will now be a super coquasi-Hopf algebra. Looking at such
an example we have the following definition and proposition:

\begin{defin}\label{qcalc} If $(H,\Phi)$ is a coquasi-Hopf algebra, we
define its first order quasi differential calculus $\Omega^1(H)$
to be $(\Omega^1,\extd)$ where $\Omega^1$ is a quasi-bimodule in
the sense
\begin{eqnarray*} a\bullet(b\bullet\tau) &=&
\Phi (a\o\tens b\o\tens \tau\o)(a\t\bullet b\t)\bullet
\tau\t\Phi^{-1} (a\thr\tens b\thr\tens\tau\thr)\\
\tau\bullet(b\bullet c) &=& \Phi (\tau\o\tens b\o\tens
c\o)(\tau\t\bullet b\t)
\bullet c\t\Phi^{-1} (\tau\thr\tens b\thr\tens c\thr)\\
a\bullet(\tau\bullet c) &=& \Phi (a\o\tens \tau\o\tens
c\o)(a\t\bullet \tau\t)\bullet c\t\Phi^{-1} (a\thr\tens
\tau\thr\tens c\thr)\end{eqnarray*} and the rest as usual (so
$\extd$ is a derivation with respect to the product $\bullet$, etc
and the calculus is bicovariant if coactions $\Delta_L,\Delta_R$
commute with $\bullet$ and intertwine $\extd$).
  \end{defin}

\begin{propos} If $H$ is an ordinary Hopf algebra and $F$ a
2-cochain then $\Omega^1(H_F)$ defined by (\ref{twicalc}) and the
same $\extd$ is a first order quasidifferential calculus on $H_F$.
\end{propos}
\proof This is more or less by construction. We define
$\Omega(H_F)$ namely as a super coquasi-Hopf algebra given by
$\Omega(H)_F$. As such it obeys the conditions in
Definition~\ref{qcalc} as these are just the lowest order part of
the assertion that $\Omega(H_F)$ is a super coquasi-Hopf algebra.
One may check that $\extd$ remains a derivation, the proof being
the same as in \cite{MaOec:twi}. \eproof

We do not go through all the steps of the Woronowicz construction
for the exterior algebra $\Omega(H)$ in detail but this may surely
be done and is straightforward in view of the category of
quasi-crossed modules being known \cite{Ma:qdou}. Then the entire
$\Omega(H_F)$ will be an example of such a general construction.

In particular, by the remarks above, all the standard
coquasitriangular quantum groups $C_\hbar(G)$ have such
bicovariant quasi-differential calculi $\Omega^1(C_\hbar(G))$ and
indeed an entire super coquasi-Hopf exterior algebra. We now make
a semiclassical analysis to lowest order of this example. What we
need from Drinfeld's theory is that
\[ F_\hbar=1\tens 1+ \hbar f+O(\hbar^2)\in U(\cg)\tens U(\cg)\]
where $\cg$ is the Lie algebra of $G$. Since it is given by
twisting from its symmetric part, the quasitriangular structure is
\[ r=f_{21}-f+r_+\]
where $r_+$ is the split Casimir or inverse Killing form in a
suitable normalisation. The antisymmetric part is $r_-=f_{21}-f$.
Moreover,
\[ \Phi=1\tens 1\tens 1+\hbar^2[[r_-,r_-]]+O(\hbar^3) \]
The following applies to these standard Poisson-Lie groups  and
any others where the quasitriangular Poisson-Lie group can be
quantized to a coquasitriangular Hopf algebra by a cochain twist,
with the induced twisted quantum differential calculus.

\begin{propos} The preconnection obtained by semiclassicaling
the quantum differential calculi obtained by cochain-twisting is
the canonical one $\hat\Xi=0$ in Theorem~\ref{cangamma}.
\end{propos}
\proof In the twisted bicovariant calculus we have
\begin{eqnarray*} a\bullet \tau&=&a\tau+\hbar(\<f,a\o\tens
\tau\bo\>a\t
\tau\bt-a\o\tau\bz\<f,a\t\tens\tau\bo\>)+O(\hbar^2)\\
\tau\bullet a&=&\tau a+\hbar(\<f,\tau\bo\tens a\o \> \tau\bt
a\t-\tau\bz a\o\<f,\tau\bo\tens
a\t\>)+O(\hbar^2)\\
{}[a,\tau]_\bullet&=&\hbar(\<f-f_{21},a\o\tens \tau\bo \>
a\t\tau\bt- a\o\tau\bz\<f-f_{21},
a\t\tens\tau\bo\>)+O(\hbar^2)\\
&=&\hbar\left( (a\ra r_-\so)(\tau\ra r_-\st)-(r_-\so\la
a)(r_-\st\la \tau)\right)+O(\hbar^2)
\end{eqnarray*} as follows. First, when we expand
$F,F^{-1}$, the $1\tens 1$ of one or other of these evaluates
against $(\Delta\tens\id)\Delta\tau$ to yield
$\Delta\tau=\Delta_L\tau+\Delta_R\tau$ and of these two terms only
will contribute as $F$ pairs only in degree zero. We recognize
$-r_-$ appearing here. Finally, we note that if $v\in\cg$ then
$\<v,a\o\>a\t=-v\la a=R_{*v}(a)$ and
$\<v,\tau\bo\>\tau\bt=-v\la\tau$ in the conventions of
Section~4.1, and similarly on the other side. Comparing, we have
exactly the formula in Proposition~\ref{cangamma}, which
corresponds to $\hat\Xi=0$. \eproof

In view of Proposition~\ref{cangamma} we also see the
nonassociativity of the exterior algebra reflected in the
curvature of the associated preconnection. Also, there is
nothing stopping one doing the above with a cochain where $\Phi$
is not central in the sense above. Then one obtains a
quasidifferential calculus on a quasi-Hopf algebra $H_F$. In this
case the infinitesimal data even for the function algebra is not a
Poisson bracket but something weaker as in \cite{Koss}.

\section{The Fedosov point of view}

The reader may wonder where the Fedosov method of deformation
quantisation \cite{Fed} fits into this picture of deforming
differential calculus that we have developed. Since the initial
data in \cite{Fed} is a symplectic structure {\em and} a
symplectic connection, it seems that its correct formulation could
be as not only quantising functions but functions and
differentials. We first point out that there are indeed flat
sections of the bundle of Weyl differential forms, and the wedge
product of such sections is flat. All multiplications are
associative, and the forms are a left and right module over the
functions on the manifold. The reader is reminded that the
functions on the manifold are replaced in the Fedosov theory by
the flat sections of the 0-form Weyl bundle, $W_D$ in the notation
of \cite{Fed}. Since these flat sections are in a natural
  1-1 correspondence with the ordinary
functions on the manifold, it seems that we have deformed
  the product on $C^\infty(M)$ itself. The problem with the
  $q$-forms for $q\ge 1$
is that there are far too many flat sections of the $q$-form Weyl
bundle, as we now see:

\begin{propos} For $q\ge 1$ there is a 1-1 correspondence
\[
\delta\delta^{-1}: \big\{ \tau\in C^\infty(W\tens\Lambda^q) :
D\tau=0\big\} \longrightarrow \{\eta\in C^\infty(W\tens\Lambda^q)
: \delta\eta=0\big\}\ .
\]
The inverse map sends $\eta$ to $\tau$, the unique solution of the
equation
\begin{eqnarray}\label{fed1}
\tau\,=\,\delta^{-1}(D+\delta)\tau\,+\,\eta\ . \end{eqnarray}
\end{propos} \proof The fact that (\ref{fed1}) has a solution can
be seen by an iterative method. Set $\tau_0=\eta$, and continue
with
\[
\tau_{n+1}\,=\,\delta^{-1}(D+\delta)\tau_n\,+\,\eta\ .
\]
Then we have
\[
\tau_{n+2}\,-\,\tau_{n+1}\,=\,
\delta^{-1}(D+\delta)(\tau_{n+1}\,-\,\tau_n)\ ,
\]
and use the fact that the operation $\delta^{-1}(D+\delta)$
increases degree to see that this iterative solution converges as
$n\to\infty$ in each degree. To show that
  (\ref{fed1}) has a unique solution, take another solution
$\tilde\tau$ and subtract to get
\[
\tilde\tau-\tau\,=\,\delta^{-1}(D+\delta)(\tilde\tau-\tau)\ .
\]
By counting degrees again, we see that $\tilde\tau=\tau$. Next we
have to show that if $\tau$ is a solution of (\ref{fed1}), then
$D\tau=0$. As $\delta^{-1}\delta^{-1}=0$, applying $\delta^{-1}$
to (\ref{fed1}) gives
\begin{eqnarray}\label{fed2}
\delta^{-1}\tau\,=\,\delta^{-1}\eta\ . \end{eqnarray} Then we
have, as $ \delta\delta^{-1} + \delta^{-1}\delta $ is the identity
on $C^\infty(W\tens\Lambda^q)$,
\begin{eqnarray}\label{fed3}
\delta^{-1}D\tau &=&\delta^{-1}(D+\delta)\tau\,-\,
\delta^{-1}\delta\tau \,=\,\tau\,-\,\eta\,-\,
\delta^{-1}\delta\tau  \cr &=& \delta\delta^{-1}\tau\,-\,\eta\,=\,
\delta\delta^{-1}\eta\,-\,\eta
\,=\,-\,\delta^{-1}\delta\eta\,=\,0\ . \end{eqnarray} Then,
remembering that $D(D\tau)=0$,
\begin{eqnarray}\label{fed4}
D\tau &=& (\delta^{-1}\delta+\delta\delta^{-1})D\tau\,=\,
\delta^{-1}\delta D\tau\,=\,\delta^{-1}(D+\delta)D\tau \ .
\end{eqnarray} By the previous degree argument, this implies that
$D\tau=0$, and concludes showing that the back map is well
defined.

To show that the composition one way round is the identity, apply
$\delta\delta^{-1}$ to (\ref{fed1}) to get
\[
\delta\delta^{-1}\tau\,=\,\delta\delta^{-1}\eta\,=\,
(\delta^{-1}\delta+\delta\delta^{-1})\eta\,=\,\eta\ .
\]

Finally we show that the composition of the maps the other way
round is also the identity. Begin with $\tau$ satisfying
$D\tau=0$. Then, as required
\[
\tau\,=\,(\delta^{-1}\delta+\delta\delta^{-1})\tau\,=\,
\delta^{-1}(D+\delta)\tau\,+\,\delta\delta^{-1}\tau\ .\quad\square
\]

\medskip
Now the classical $q$-forms do give rise to flat Weyl $q$-forms,
as $C^\infty(M)\tens \Lambda^q$ is contained in the kernel of
$\delta$, so we can use (\ref{fed1}) to find a corresponding flat
section. However the kernel of $\delta$ in
$C^\infty(W\tens\Lambda^q)$ is much larger
  than the classical $q$-forms. If we take the wedge product of two flat
sections coming from classical $q$-forms, then we get a flat
section, but not necessarily one coming from a classical $q$-form.

We are then faced with a choice, either to stick with the
associative framework given by flat Weyl forms and sacrifice
correspondence with the classical forms, or to try to maintain the
correspondence but accept that associativity will only be
approximately true.

Remember that a function $f\in C^{\infty}(M)$ has a corresponding
quantisation beginning $f+y^{k}\,\nabla_{k}f+\frac12
y^{k}y^{j}\,\nabla_{k}\nabla_{j}f+\dots$, so if we had a 1-form
also expanded in the form $\eta+y^{j}\eta_{[j]}+\dots$, we would
get
\[
[f+y^{i}\,\nabla_{i}f+\dots,\eta+y^{i}\eta_{[i]}+\dots]\,=\,
-\,i\,\hbar\,\omega^{kj}\,\eta_{[j]}\,\nabla_{k}f\,+\,\dots\ ,
\]
so (using $-i\hbar$ instead of $\hbar$ to fit the usual notation
of the Fedosov theory), we would have to have the quantisation of
the 1-form beginning $\eta+y^{k}\,\nabla_{k}\eta+\dots$.  However
a calculation of the leading order expansion of the flat 1-form
whose lowest order part is $\eta$ gives
$\eta_{n}\,+\,\frac12\,y^{k}(\eta_{n;k}-\eta_{k;n})+\dots$. We
conclude that taking flat 1-forms is really not the right thing to
do.  Also note that the super-commutator of two 1-forms is, to
leading order,
\[
[\eta+y^{j}\,\nabla_{j}\eta+\dots,\xi+y^{k}\,\nabla_{k}\xi+\dots]\,=\,
-\,i\,\hbar\,\omega^{jk}\,\nabla_{j}\eta\wedge\nabla_{k}\xi\ ,
\]
also as previously described.


\end{document}